\documentclass{amsart}
\usepackage{amscd,amssymb}
\usepackage{righttag}
\usepackage{epsfig}
\usepackage{pb-diagram,lamsarrow,pb-lams}

\newcommand{\Aut}       {\operatorname{Aut}}
\newcommand{\Level}     {\operatorname{Level}}
\newcommand{\Out}       {\operatorname{Out}}
\newcommand{\Hom}       {\operatorname{Hom}}
\newcommand{\Sets}      {\operatorname{Sets}}

\newcommand{\spec}      {\operatorname{spec}}
\newcommand{\spf}       {\operatorname{spf}}
\newcommand{\image}     {\operatorname{image}}
\newcommand{\mor}       {\operatorname{mor}}
\newcommand{\obj}       {\operatorname{obj}}
\newcommand{\tr}        {\operatorname{tr}}

\newcommand{\C}         {{\mathbb{C}}}
\newcommand{\R}         {{\mathbb{R}}}
\newcommand{\Z}         {{\mathbb{Z}}}              
\newcommand{\Q}         {{\mathbb{Q}}}
\newcommand{\FF}        {{\mathbb{F}}}
\newcommand{\Qp}        {{\mathbb{Q}_p}}          
\newcommand{\Zp}        {{\mathbb{Z}_p}}          

\newcommand{\al}        {\alpha}
\newcommand{\bt}        {\beta} 
\newcommand{\gm}        {\gamma}
\newcommand{\dl}        {\delta}
\newcommand{\ep}        {\epsilon}
\newcommand{\tht}       {\theta}
\newcommand{\lm}        {\lambda}
\newcommand{\sg}        {\sigma}
\newcommand{\zt}        {\zeta}

\newcommand{\Lm}        {\Lambda}
\newcommand{\Sgi}       {\Sigma^\infty}
\newcommand{\Sg}        {\Sigma}
\newcommand{\Om}        {\Omega}

\newcommand{\CB}        {{\mathcal{B}}}
\newcommand{\CBb}       {{\overline{\mathcal{B}}}}
\newcommand{\CBf}       {{\mathcal{B}_f}}
\newcommand{\CBbf}      {{\overline{\mathcal{B}}_f}}
\newcommand{\CC}        {{\mathcal{C}}}
\newcommand{\CF}        {{\mathcal{F}}}
\newcommand{\CG}        {{\mathcal{G}}}
\newcommand{\CGb}       {{\overline{\mathcal{G}}}}
\newcommand{\CGf}       {{\mathcal{G}_f}}
\newcommand{\CGbf}      {{\overline{\mathcal{G}}_f}}
\newcommand{\DD}        {{\mathcal{D}}}
\newcommand{\CK}        {{\mathcal{K}}}
\newcommand{\CM}        {{\mathcal{M}}}
\newcommand{\CS}        {{\mathcal{S}}}
\newcommand{\ot}        {\otimes}
\newcommand{\Smash}     {\wedge}
\newcommand{\Wedge}     {\vee}
\newcommand{\bigWedge}  {\bigvee}
\newcommand{\tm}        {\times}
\newcommand{\xra}       {\xrightarrow}
\newcommand{\xla}       {\xleftarrow}
\newcommand{\st}        {\;|\;}
\newcommand{\GG}        {{\mathbb{G}}}
\newcommand{\mxi}       {\mathfrak{m}}
\newcommand{\tq}        {\tilde{q}}
\newcommand{\tv}        {\tilde{v}}
\newcommand{\eqdef}     {:=}
\newcommand{\sse}       {\subseteq}
\newcommand{\gb}        {\overline{g}}

\newcommand{\Dl}        {\Delta}
\newcommand{\Gm}        {\Gamma}
\newcommand{\kp}        {\kappa}
\newcommand{\ttt}       {\tilde{\tau}}
\renewcommand{\O}       {{\mathcal{O}}}
\newcommand{\cpi}       {{\mathbb{C}P^\infty}}
\newcommand{\dps}[1]    {\langle p^#1\rangle}
\newcommand{\fps}[2]    {{#1 [\![ #2 ]\!]}} 
\newcommand{\om}        {\omega}
\newcommand{\res}       {\operatorname{res}}
\newcommand{\tE}        {\widetilde{E}}
\newcommand{\hS}        {\widehat{S}}
\newcommand{\hSmash}    {\widehat{\wedge}}

\newcommand{\colim}  {\operatornamewithlimits{\underset{\longrightarrow}{lim}}}
\newcommand{\invlim} {\operatornamewithlimits{\underset{\longleftarrow}{lim}}}

\renewcommand{\:}{\colon}

\newtheorem{theorem}{Theorem}[section]

\newtheorem{lemma}[theorem]{Lemma}
\newtheorem{proposition}[theorem]{Proposition}
\newtheorem{corollary}[theorem]{Corollary}
\newtheorem{scholium}[theorem]{Scholium}
\theoremstyle{definition}
\newtheorem{remark}[theorem]{Remark}
\newtheorem{definition}[theorem]{Definition}
\newtheorem{example}[theorem]{Example}
\newtheorem{convention}[theorem]{Convention}

\newenvironment{diag}{
 \renewcommand{\typeout}[1]{}
 \begin{displaymath}
 \begin{diagram}}{
 \end{diagram}
 \end{displaymath}} 

\begin{document}
\title{{$K(n)$}-local duality for finite groups and groupoids}
\author{N.~P.~Strickland}
\date{\today}
\bibliographystyle{abbrv}

\maketitle 

\section{Introduction}

The starting point of the investigations described here was our
discovery of a natural inner product on the ring $K(n)^*BG$, the
$n$'th Morava $K$-theory of the classifying space of a finite group
$G$.  If $n=1$ and $G$ is a $p$-group then $K(1)^*BG$ is essentially
the same as $R(G)/p$ (where $R(G)$ is the complex representation ring
of $G$) and our inner product is just 
$(V,W)=\dim_\C(V\ot W)^G\pmod{p}$.  This is closely related to the
classical inner product on $R(G)$, given by
\[ \langle V,W\rangle = \dim_\C\Hom_G(V,W) = (V,W^*). \]
(For more general groups $G$, there is still a relationship with the
classical product but it is not too close; see
Section~\ref{sec-warnings} for some pitfalls.)  

It turns out to be useful to work with an inner product on the
spectrum $LG:=L_{K(n)}\Sgi BG_+$ and then deduce consequences in
Morava $K$-theory (and other generalised cohomology theories) by
functorality.  As background to this, in Section~\ref{sec-inner} we
recall some results about inner products on objects in arbitrary
compact closed categories.  Moreover, to elucidate the relationship
between the inner product and the ring structure on $K(n)^*BG$, it is
helpful to recall some facts about Frobenius algebras and their
relationship with topological quantum field theories, which we do in
Sections~\ref{sec-frobenius} and~\ref{sec-trace-form}.  In
Section~\ref{sec-manifolds} we give a version of Poincar\'{e}-Atiyah
duality for manifolds which illustrates these ideas nicely, and which
has striking formal similarities with our later treatment of $LG$;
indeed, one could probably set up a unifying categorical framework.
We have also found that many aspects of our theory (for example
homotopy pullbacks and free loop spaces) can be discussed more cleanly
in terms of groupoids rather than groups.  This is also convenient for
a number of applications and calculations.  Because of this, we give a
fairly detailed treatment of the homotopy theory of groupoids in
Section~\ref{sec-groupoids}.  In Section~\ref{sec-transfers} we
discuss transfers for coverings-up-to-homotopy, as outlined
in~\cite[Remark 3.1]{sttu:rme}.  In Section~\ref{sec-Kn-local} we turn
to the spectra $LG$.  In~\cite{host:mkl} we used the Greenlees-May
theory of generalised Tate spectra to exhibit an equivalence 
$LG\simeq F(LG,L_{K(n)}S^0)$.  After comparing some definitions and
feeding this into our machinery, we find that $LG$ has a natural
structure as a Frobenius object in the $K(n)$-local stable category,
whenever $G$ is a finite groupoid.  As part of the construction we
define $K(n)$-local transfer maps for arbitrary homomorphisms of
finite groups, or functors of finite groupoids; these reduce to
classical transfers when the homomorphisms or functors are injective
or faithful.  In Section~\ref{sec-inner-cohomology} we deduce various
consequences for the generalised cohomology of $BG$; in the case where
$G$ is a finite Abelian group, we can be quite explicit.  In
Section~\ref{sec-characters}, we deduce some further consequences in
terms of the Hopkins-Kuhn-Ravenel generalised character
theory~\cite{hokura:ggc}, which gives a complete description of 
$\Q\ot E^0BG$ for suitable cohomology theories $E$.  Finally, in
Section~\ref{sec-warnings} we alert the reader to some possible
pitfalls that can arise from overoptimism about the analogy with
classical representation theory.

\section{Inner products}
\label{sec-inner}

Let $\CC$ be an additive compact closed category, in other words an
additive closed symmetric monoidal category in which every object is
dualisable.  We write $X\Smash Y$ for the symmetric monoidal product,
and $S$ for the unit object.  We also write $F(Y,Z)$ for the function
objects, so that $\CC(X,F(Y,Z))\simeq\CC(X\Smash Y,Z)$.  Finally, we
write $DX=F(X,S)$, so that $D^2X=X$ and $F(X,Y)=DX\Smash Y$.

\begin{definition}
 An \emph{inner product} on an object $X\in\CC$ is a map
 $b\:X\Smash X\xra{}S$ such that 
 \begin{enumerate}
 \item $b$ is symmetric in the sense that $b\circ\tau=b$, where
  $\tau\:X\Smash X\xra{}X\Smash X$ is the twist map; and
  \item the adjoint map $b^\#\:X\xra{}DX$ is an isomorphism.
 \end{enumerate}
\end{definition}

\begin{example}
 We could take $\CC$ to be the category of finitely generated
 projective modules over a commutative ring $R$, with the usual closed
 symmetric monoidal structure so that $M\Smash N=M\ot_RN$ and
 $F(M,N)=\Hom_R(M,N)$ and $DM=M^*=\Hom_R(M,R)$.  An inner product on
 $M$ is then a symmetric $R$-bilinear pairing $M\tm M\xra{}R$ that
 induces an isomorphism $M\simeq M^*$.  If $R$ is a field then this
 just says that the pairing is nondegenerate.  Note that we have no
 positivity condition.
\end{example}

\begin{remark}\label{rem-duality-equational}
 We see from~\cite[Theorem III.1.6]{lemast:esh} that a symmetric map
 $b\:X\Smash X\xra{}S$ is an inner product iff it is a duality of $X$
 with itself in the sense discussed there, iff there is a map
 $c\:S\xra{}X\Smash X$ such that the following diagrams commute:
 \begin{diag}
  \node{X} \arrow{e,t}{1\Smash c} \arrow{se,b}{1}
  \node{X\Smash X\Smash X} \arrow{s,r}{b\Smash 1}
  \node[2]{X\Smash X\Smash X} \arrow{s,l}{1\Smash b}
  \node{X} \arrow{w,t}{c\Smash 1} \arrow{sw,b}{1} \\
  \node[2]{X} \node[2]{X}
 \end{diag}
 Moreover, if $b$ is an inner product then there is a unique map $c$
 as above, and it is symmetric; in fact it is also the unique
 symmetric map making the left hand diagram commute.
\end{remark}
\begin{remark}
 The commutativity of the above diagrams can be expressed in terms of
 Penrose diagrams~\cite{jost:gtci} as follows:
 \[
 \setlength{\unitlength}{0.0005in}%
 \begin{picture}(6624,785)(589,-463)
  \put(1801,239){\circle*{100}}
  \put(1201,-361){\circle*{100}}
  \put(6601,-361){\circle*{100}}
  \put(6001,239){\circle*{100}}
  \put(601,239){\line( 1, 0){1200}}
  \put(1801,239){\line(-1,-1){600}}
  \put(1201,-361){\line( 1, 0){1200}}
  \put(3301,-61){\line( 1, 0){1200}}
  \put(5401,-361){\line( 1, 0){1200}}
  \put(6601,-361){\line(-1, 1){600}}
  \put(6001,239){\line( 1, 0){1200}}
  \put(2776,-136){\makebox(0,0){$=$}}
  \put(4876,-136){\makebox(0,0){$=$}}
  \put(2026,164){\makebox(0,0) {$b$}}
  \put(826,-436){\makebox(0,0) {$c$}}
  \put(5626,164){\makebox(0,0) {$c$}}
  \put(6826,-436){\makebox(0,0){$b$}}
  \put(3826,-361){\makebox(0,0){$1$}}
 \end{picture}
 \]
 Similarly, the symmetry of $b$ and $c$ gives the following equations:
 \[
 \setlength{\unitlength}{0.0005in}%
 \begin{picture}(8124,1224)(589,-973)
  \put(1201,-361){\circle*{100}}
  \put(2401,-361){\circle*{100}}
  \put(3901,-361){\circle*{100}}
  \put(5401,-361){\circle*{100}}
  \put(6901,-361){\circle*{100}}
  \put(8101,-361){\circle*{100}}
  \put(601,239){\line( 1,-1){1200}}
  \put(1801,-961){\line( 1, 1){600}}
  \put(2401,-361){\line(-1, 1){600}}
  \put(1801,239){\line(-1,-1){1200}}
  \put(3301,239){\line( 1,-1){600}}
  \put(3901,-361){\line(-1,-1){600}}
  \put(6001,239){\line(-1,-1){600}}
  \put(5401,-361){\line( 1,-1){600}}
  \put(8701,-961){\line(-1, 1){1200}}
  \put(7501,239){\line(-1,-1){600}}
  \put(6901,-361){\line( 1,-1){600}}
  \put(7501,-961){\line( 1, 1){1200}}
  \put(2776,-436){\makebox(0,0){$=$}}
  \put(6376,-436){\makebox(0,0){$=$}}
  \put(4650,-361){\makebox(0,0){and}}
  \put( 826,-436){\makebox(0,0){$\tau$}}
  \put(2026,-436){\makebox(0,0){$b$}}
  \put(3526,-436){\makebox(0,0){$b$}}
  \put(5701,-436){\makebox(0,0){$c$}}
  \put(7201,-436){\makebox(0,0){$c$}}
  \put(8401,-436){\makebox(0,0){$\tau$}}
 \end{picture}
 \]
\end{remark}

\begin{definition}
 If $X$ and $Y$ are equipped with inner products and $f\:X\xra{}Y$
 then we write $f^t\:Y\xra{}X$ for the unique map making the following
 diagram commute:
 \begin{diag}
  \node{Y}  \arrow{s,l}{b_Y^\#} \arrow{e,t}{f^t} 
  \node{X}  \arrow{s,r}{b_X^\#}                  \\
  \node{DY} \arrow{e,b}{Df}     \node{DX.}
 \end{diag}
 This can also be characterised by the equation
 \[ b_Y\circ(f\Smash 1) = b_X\circ(1\Smash f^t) \: X\Smash Y\xra{}S \]
 or equivalently, the following equality between Penrose diagrams.
 \[
 \setlength{\unitlength}{0.0005in}
 \begin{picture}(5108,1410)(976,-1138)
  \put(2101,-61){\circle*{100}}
  \put(3001,-361){\circle*{100}}
  \put(6001,-361){\circle*{100}}
  \put(5101,-661){\circle*{100}}
  \put(1201,239){\line( 3,-1){1800}}
  \put(3001,-361){\line(-3,-1){1800}}
  \put(4201,239){\line( 3,-1){1800}}
  \put(6001,-361){\line(-3,-1){1800}}
  \put(3526,-436) {\makebox(0,0){$=$}}
  \put(2851,-736) {\makebox(0,0){$b_Y$}}
  \put(2101,164)  {\makebox(0,0){$f$}}
  \put(5101,-1111){\makebox(0,0){$f^t$}}
  \put(5851,-736) {\makebox(0,0){$b_X$}}
  \put(3976,164)  {\makebox(0,0){$X$}}
  \put(3976,-1036){\makebox(0,0){$Y$}}
  \put(976,164)   {\makebox(0,0){$X$}}
  \put(976,-1036) {\makebox(0,0){$Y$}}
 \end{picture}
 \]
 It is clear that $f^{tt}=f$ and that $1^t=1$ and $(gf)^t=f^tg^t$
 whenever this makes sense.  We call $f^t$ the \emph{transpose} of
 $f$.
\end{definition}

\begin{remark}
 Suppose that $X$ and $Y$ have inner products $b_X$ and $b_Y$.  We
 then define
 \[ b_{X\Smash Y} =
     (X\Smash Y\Smash X\Smash Y \xra{1\Smash\tau\Smash 1}
     X\Smash X\Smash Y\Smash Y \xra{b_X\Smash b_Y} S).
 \]
 It is easy to check that this is an inner product on $X\Smash Y$.
 Similarly, if $\CC$ is an additive category (with direct sums written
 as $X\Wedge Y$) and $\Smash$ is bilinear then there is an obvious way
 to put an inner product on $X\Wedge Y$.  By abuse of language, we
 call these inner products $b_X\Smash b_Y$ and $b_X\Wedge b_Y$.  If we
 use these inner products, we find that $(f\Smash g)^t=f^t\Smash g^t$
 and $(f\Wedge g)^t=f^t\Wedge g^t$.
\end{remark}

\section{Frobenius objects}
\label{sec-frobenius}

\begin{definition}
 Let $\CC$ be a symmetric monoidal category.  A \emph{Frobenius
   object} in $\CC$ is an object $A\in\CC$ equipped with maps
 $S\xra{\eta}A\xla{\mu}A\Smash A$ and $S\xla{\ep}A\xra{\psi}A\Smash A$
 such that
 \begin{itemize}
  \item[(a)] $(A,\eta,\mu)$ is a commutative and associative ring
   object. 
  \item[(b)] $(A,\ep,\psi)$ is a commutative and associative coring
   object. 
  \item[(c)] (The ``interchange axiom'') The following diagram
   commutes:
   \begin{diag}
    \node{A\Smash A} \arrow{e,t}{\mu} \arrow{s,l}{\psi\Smash 1}
    \node{A}                          \arrow{s,r}{\psi} \\
    \node{A\Smash A\Smash A} \arrow{e,b}{1\Smash\mu}
    \node{A\Smash A.}
   \end{diag}
 \end{itemize}
\end{definition}
The point for us will be that for any finite groupoid $G$, the
spectrum $LG:=L_{K(n)}\Sgi BG_+$ has a natural structure as a
Frobenius object in the $K(n)$-local stable category
(Theorem~\ref{thm-groupoids-frobenius}). 

The last axiom can be restated as the following equality of Penrose
diagrams: 
\[
 \setlength{\unitlength}{0.0005in}%
 \begin{picture}(4824,1366)(589,-1044)
 \put(1201,-361){\circle*{100}}
 \put(1280,-560){\makebox(0,0){$\mu$}}
 \put(1801,-361){\circle*{100}}
 \put(1720,-560){\makebox(0,0){$\psi$}}
 \put(3901,239){\circle*{100}}
 \put(3850, 40){\makebox(0,0){$\psi$}}
 \put(5101,-961){\circle*{100}}
 \put(5150,-760){\makebox(0,0){$\mu$}}
 \put(2800,-361){\makebox(0,0){$=$}}
 \put(601,239){\line( 1,-1){600}}
 \put(1201,-361){\line(-1,-1){600}}
 \put(1201,-361){\line( 1, 0){600}}
 \put(2401,239){\line(-1,-1){600}}
 \put(1801,-361){\line( 1,-1){600}}
 \put(3601,239){\line( 1, 0){1800}}
 \put(3601,-961){\line( 1, 0){1800}}
 \put(3901,239){\line( 1,-1){1200}}
 \end{picture}
\]
\begin{remark}\label{rem-axioms-self-dual}
 If $(A,\eta,\ep,\mu,\psi)$ is a Frobenius object in $\CC$ then
 $(A,\ep,\eta,\psi,\mu)$ is evidently a Frobenius object in
 $\CC^{\text{op}}$. 
\end{remark}
\begin{convention}\label{conv-penrose}
 For the rest of this paper, we use the following conventions for
 Penrose diagrams.  Unless otherwise specified, each diagram will
 involve only a single object $A$, for which some subset of the maps
 $\mu$, $\psi$, $\eta$, $\ep$ will have been defined.  We also
 automatically have a twist map $\tau\:A\Smash A\xra{}A\Smash A$.
 \begin{itemize}
  \item Any unlabelled node with two lines in and one line out is
   implicitly labelled with $\mu$.
  \item A node with one line in and no lines out is implicitly
   labelled $\psi$.
  \item A node with no lines in and one line out is implicitly
   labelled $\eta$.
  \item A node with one line in and no lines out is implicitly
   labelled $\ep$.
  \item A node with two lines in and two lines out is implicitly
   labelled $\tau$.
 \end{itemize}
\end{convention}

Another interesting point of view is that Frobenius objects are
equivalent to topological quantum field theories (TQFT's).  In more
detail, let $\CS$ be the $1+1$-dimensional cobordism category, whose
objects are closed $1$-manifolds and whose morphisms are cobordisms.
Some care is needed to set the details up properly: a good account
is~\cite{ab:tdt}, although apparently the results involved were ``folk
theorems'' long before this.  The category $\CS$ has a symmetric
monoidal structure given by disjoint unions.  The circle $S^1$ is a
Frobenius object in $\CS$: the maps $\eta$ and $\ep$ are the disc
$D^2$ regarded as a morphism $\emptyset\xra{}S^1$ and
$S^1\xra{}\emptyset$ respectively, and the maps $\mu$ and $\psi$ are
the ``pair of pants'' regarded as a morphism $S^1\amalg S^1\xra{}S^1$
and $S^1\xra{}S^1\amalg S^1$ respectively.  It follows easily
from~\cite[Proposition 12]{ab:tdt} that this is a universal example of
a symmetric monoidal category equipped with a Frobenius object.  For
further analysis of the category $\CS$, see~\cite{ca:cc,ti:csi}.

\begin{remark}
 Using the Frobenius structure on $S^1$, a Penrose diagram as in
 Convention~\ref{conv-penrose} gives rise to a morphism in $\CS$.
 This has the following appealing geometric interpretation.  We first
 perform the replacement 
 \[
 \setlength{\unitlength}{0.0005in}%
 \begin{picture}(6924,1824)(589,-1573)
  \put(1501,-661){\circle*{100}}
  \put(601,239){\line( 1,-1){1800}}
  \put(601,-1561){\line( 1, 1){1800}}
  \put(3301,-361){\line( 0,-1){600}}
  \put(3301,-661){\line( 1, 0){1500}}
  \put(4501,-361){\line( 1,-1){300}}
  \put(4801,-661){\line(-1,-1){300}}
  \put(5701,239){\line( 1,-1){1800}}
  \put(5701,-1561){\line( 1, 1){675}}
  \put(7501,239){\line(-1,-1){750}}
 \end{picture}
 \]
 (It makes no real difference whether we introduce an under crossing
 or an over crossing.)  This converts the Penrose diagram to a graph
 embedded in $[0,1]\tm\R^2$.  The boundary of a regular neighbourhood
 of this graph is a surface $\Sg$ which we can think of as a cobordism
 between $\Sg\cap(\{0\}\tm\R^2)$ and $\Sg\cap(\{1\}\tm\R^2)$ and thus
 as a morphism in $\CS$.  For example, the Penrose diagram
 \[
 \setlength{\unitlength}{0.0005in}%
 \begin{picture}(3624,4224)(589,-3373)
  \put(1201,-1561){\circle*{100}}
  \put(1501,-1561){\circle*{100}}
  \put(1801,-1561){\circle*{100}}
  \put(3001,-1561){\circle*{100}}
  \put(3601,-1561){\circle*{100}}
  \put(2101,-1861){\circle*{100}}
  \put(2701,-1261){\circle*{100}}
  \put(1201,-2761){\circle*{100}}
  \put(2401,-2761){\circle*{100}}
  \put(601,-3361){\line( 1, 1){3600}}
  \put(4201,-961){\line(-1,-1){600}}
  \put(3601,-1561){\line( 1,-1){600}}
  \put(3001,-1561){\line( 1, 0){600}}
  \put(1801,-1561){\line( 1, 1){600}}
  \put(2401,-961){\line( 1,-1){600}}
  \put(3001,-1561){\line(-1,-1){600}}
  \put(2401,-2161){\line(-1, 1){600}}
  \put(601,-2161){\line( 1, 1){600}}
  \put(1201,-1561){\line(-1, 1){600}}
  \put(1201,-1561){\line( 1, 0){600}}
  \put(601,239){\line( 1,-2){900}}
  \put(1201,-2761){\line( 1, 0){1200}}
 \end{picture}
 \]
 becomes the following cobordism:
 \[ \epsfbox{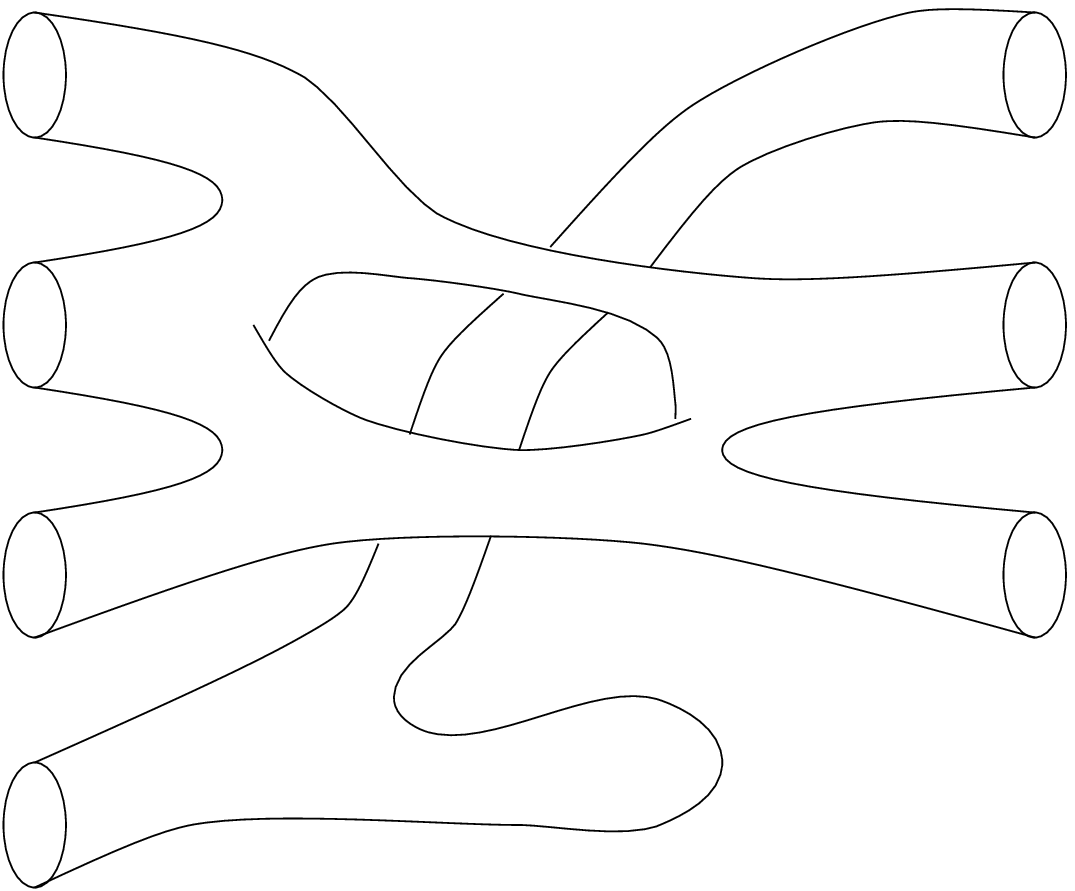} \] 
\end{remark}

\begin{definition}
 Let $\CC$ be a compact closed category, and let $A$ be an object of
 $\CC$ equipped with a commutative and associative product
 $\mu\:A\Smash A\xra{}A$.  (We do not assume that there is a unit.)  A
 \emph{Frobenius form} on $A$ is a map $\ep\:A\xra{}S$ such that the
 map $b=\ep\mu$ is an inner product.
\end{definition}
\begin{example}
 The most familiar example in topology is that if $M$ is a closed
 manifold with fundamental class $[M]\in H_*(M;\FF_2)$ then the
 equation $\ep(u)=\langle u,[M]\rangle$ is a Frobenius form on
 $H^*(M;\FF_2)$ (regarded as an ungraded module over $\FF_2$).  This
 can of course be generalised to other coefficients at the price of a
 few words about orientations and gradings.  For a geometrised version
 of this, see Section~\ref{sec-manifolds}.
\end{example}
\begin{example}
 Another elementary example is to let $k$ be a field and $G$ a finite
 Abelian group.  We can then define a map $\ep\:k[G]\xra{}k$ sending
 $[1]$ to $1$ and $[g]$ to $0$ for $g\neq 1$.  This is easily seen to
 be a Frobenius form.
\end{example}
\begin{example}
 Let $\CC$ be the category of finitely generated free Abelian groups.
 Let $G$ be a finite group, let $R=R(G)$ be its complex representation
 ring, and define $\ep\:R\xra{}\Z$ by $\ep[W]=\dim_\C W^G$.  It is
 easy to see that this is a Frobenius form, and that the associated
 inner product is $([U],[W])=\dim_\C(U\ot W)^G$ as considered in the
 introduction.  To generalise this to finite groupoids, let $V$ be the
 category of finite dimensional complex vector spaces.  A
 \emph{representation} of $G$ means a functor $G\xra{}V$.  The set
 $R_+(G)=\pi_0[G,V]$ of isomorphism classes of representations has a
 natural structure as a semiring, and we let $R(G)$ denote its group
 completion.  If $W$ is a representation then we write
 $W^G:=\invlim_GW\in V$ and $t[W]=\dim_\C W^G$ as before.  One can
 easily deduce from the classical case that this is a Frobenius form
 on $R(G)$.
\end{example}

\begin{lemma}
 If $(A,\eta,\ep,\mu,\psi)$ is a Frobenius object in a compact closed
 category $\CC$, then $\ep$ is a Frobenius form.
\end{lemma}
\begin{proof}
 Put $b=\ep\mu\:A\Smash A\xra{}S$; we need to show that this is an
 inner product.  Put $c=\psi\eta\:S\xra{}A\Smash A$; it will suffice
 to check the identities in Remark~\ref{rem-duality-equational}.  The
 symmetry conditions are clear, so we just need the two compatibility
 conditions for $b$ and $c$.  One of them is proved as follows:
 \[
 \setlength{\unitlength}{0.0005in}%
 \begin{picture}(8424,1335)(1489,-1063)
  \put(4501,  -61){\circle*{100}}
  \put(4801,  -61){\circle*{100}}
  \put(5401, -661){\circle*{100}}
  \put(5701, -661){\circle*{100}}
  \put(6901,  -61){\circle*{100}}
  \put(7201, -361){\circle*{100}}
  \put(7801, -361){\circle*{100}}
  \put(8101, -661){\circle*{100}}
  \put(2101,  -61){\circle*{100}}
  \put(2701, -661){\circle*{100}}
  \put(4501,  -61){\line( 1, 0){1500}}
  \put(4201, -661){\line( 1, 0){1500}}
  \put(4801,  -61){\line( 1,-1){ 600}}
  \put(6901, -661){\line( 1, 1){ 300}}
  \put(7201, -361){\line(-1, 1){ 300}}
  \put(7201, -361){\line( 1, 0){ 600}}
  \put(7801, -361){\line( 1, 1){ 300}}
  \put(7801, -361){\line( 1,-1){ 300}}
  \put(9001, -361){\line( 1, 0){ 900}}
  \put(3301,  -61){\line(-1, 0){1200}}
  \put(2101,  -61){\line( 1,-1){ 600}}
  \put(2701, -661){\line(-1, 0){1200}}
  \put(3676, -436){\makebox(0,0){$=$}}
  \put(6376, -436){\makebox(0,0){$=$}}
  \put(8476, -436){\makebox(0,0){$=$}}
  \put(9376, -661){\makebox(0,0){$1$}}
  \put(1901,  -61){\makebox(0,0){$c$}}
  \put(2901, -661){\makebox(0,0){$b$}}
 \end{picture}
 \]
 The first equation is just the definition of $b$ and $c$, the second
 is the interchange axiom, and the third uses the (co)unit properties
 of $\eta$ and $\ep$.  The other compatibility condition follows
 because $b$ and $c$ are symmetric.
\end{proof}

We now prove a converse to the above result.
\begin{proposition}\label{prop-psi}
 Let $\CC$ be a compact closed category, and let $A$ be an object
 equipped with a commutative and associative product $\mu$ and a
 Frobenius form $\ep$.  Then there are unique maps $\eta,\psi$ making
 $(A,\eta,\ep,\mu,\psi)$ into a Frobenius object.
\end{proposition}
\begin{proof}
 By hypothesis $b=\ep\mu$ is an inner product on $A$, and trivially
 the canonical isomorphism $S\Smash S=S$ is an inner product on $S$.
 We can thus define $\eta\eqdef\ep^t\:S\xra{}A$, so $\eta$ is the
 unique map such that $b\circ(1\Smash\eta)=\ep$, or in other words the
 unique map giving the following equality of Penrose diagrams:
 \[
 \setlength{\unitlength}{0.0005in}%
 \begin{picture}(3695,766)(589,-444)
  \put(901,-361){\circle*{100}}
  \put(1501,239){\circle*{100}}
  \put(2101,239){\circle*{100}}
  \put(4201,239){\circle*{100}}
  \put( 601, 239){\line( 1, 0){1500}}
  \put( 901,-361){\line( 1, 1){600}}
  \put(3301, 239){\line( 1, 0){900}}
  \put(2626,-136){\makebox(0,0){$=$}}
  \put( 676,-436){\makebox(0,0){$\eta$}}
 \end{picture}
 \]
 We claim that $\eta$ is a unit for $\mu$, or in other words that we
 have the following equality:
 \[
 \setlength{\unitlength}{0.0005in}%
 \begin{picture}(3624,776)(589,-454)
  \put( 901,-361){\circle*{100}}
  \put(1501, 239){\circle*{100}}
  \put( 601, 239){\line( 1, 0){1500}}
  \put( 901,-361){\line( 1, 1){600}}
  \put(3301, 239){\line( 1, 0){900}}
  \put(2626,-136){\makebox(0,0){$=$}}
  \put(3676, -61){\makebox(0,0){$1$}}
 \end{picture}
 \]
 To prove this, we observe that for any two maps $f,g\:B\xra{}A$ we
 have $f=g$ if and only if
 \[ b\circ(1\Smash f)=b\circ(1\Smash g)\:A\Smash B\xra{}S. \]
 In view of this, the claim follows from the following diagram, in
 which the first equality comes from the associativity of $\mu$ and
 the second from the defining property of $\eta$.
 \[
 \setlength{\unitlength}{0.0005in}%
 \begin{picture}(6695,1515)(589,-1054)
  \put( 901, -961){\circle*{100}}
  \put(1501, -361){\circle*{100}}
  \put(2401,  239){\circle*{100}}
  \put(4201,  239){\circle*{100}}
  \put(3901, -961){\circle*{100}}
  \put(6901,  239){\circle*{100}}
  \put(2701,  239){\circle*{100}}
  \put(5401,  239){\circle*{100}}
  \put(7201,  239){\circle*{100}}
  \put(5101,  239){\circle*{100}}
  \put( 601,  239){\line( 1, 0){2100}}
  \put( 601, -361){\line( 1, 0){1200}}
  \put(1801, -361){\line( 1, 1){600}}
  \put( 901, -961){\line( 1, 1){600}}
  \put(3601,  239){\line( 1, 0){1800}}
  \put(3601, -361){\line( 1, 1){600}}
  \put(3901, -961){\line( 1, 1){1200}}
  \put(6301,  239){\line( 1, 0){900}}
  \put(6301, -361){\line( 1, 1){600}}
  \put(2926, -436){\makebox(0,0){$=$}}
  \put(5626, -436){\makebox(0,0){$=$}}
 \end{picture}
 \]
 We next equip $A\Smash A$ with the inner product $b\Smash b$ and
 define $\psi\eqdef\mu^t$.  As $(A,\eta,\mu)$ is a commutative and
 associative monoid object, it is easy to deduce that
 $(A,\ep,\psi)=(A,\eta^t,\mu^t)$ is a commutative and associative
 comonoid object.  Thus, to prove that $A$ is a Frobenius object, we
 need only check the interchange axiom.

 It follows directly from the definition that $\psi$ is the unique map
 giving the following equality:
 \[
 \setlength{\unitlength}{0.0005in}%
 \begin{picture}(4295,1366)(589,-1044)
  \put(1201,-361){\circle*{100}}
  \put(1801,-961){\circle*{100}}
  \put(1801, 239){\circle*{100}}
  \put(2401,-961){\circle*{100}}
  \put(2401, 239){\circle*{100}}
  \put(4201,-361){\circle*{100}}
  \put(4401,-361){\circle*{100}}
  \put(4801,-361){\circle*{100}}
  \put(1120,-560){\makebox(0,0){$\psi$}}
  \put(3000,-361){\makebox(0,0){$=$}}
  \put( 601,-361){\line( 1, 0){600}}
  \put(1201,-361){\line( 1, 1){600}}
  \put(1201,-361){\line( 1,-1){600}}
  \put( 601, 239){\line( 1, 0){1800}}
  \put( 601,-961){\line( 1, 0){1800}}
  \put(3601, 239){\line( 1, 0){200}}
  \put(3801, 239){\line( 1,-1){600}}
  \put(4201,-361){\line(-1,-1){600}}
  \put(3601,-361){\line( 1, 0){1200}}
 \end{picture}
 \]
 Using the perfectness of $b$, we see that two maps
 $f,g\:B\xra{}A\Smash A$ are equal if and only if we have
 \[ (b\Smash b)(1\Smash f\Smash 1) = (b\Smash b)(1\Smash g\Smash 1)
     \: A\Smash B\Smash A \xra{} A\Smash A\Smash A\Smash A. 
 \]
 In view of this, the interchange axiom is equivalent to the following
 equation:
 \[
 \setlength{\unitlength}{0.0005in}%
 \begin{picture}(4595,1366)(589,-1044)
 \put( 901,-361){\circle*{100}}
 \put(1201,-361){\circle*{100}}
 \put(1801, 239){\circle*{100}}
 \put(2101, 239){\circle*{100}}
 \put(1801,-961){\circle*{100}}
 \put(2101,-961){\circle*{100}}
 \put(3601, -61){\circle*{100}}
 \put(4201, -661){\circle*{100}}
 \put(4801,-961){\circle*{100}}
 \put(5101,-961){\circle*{100}}
 \put(5101, 239){\circle*{100}}
 \put(4801, 239){\circle*{100}}
 \put(2800,-361){\makebox(0,0){$=$}}
 \put( 901,-361){\line( 1, 0){300}}
 \put(1801, 239){\line(-1,-1){600}}
 \put(1201,-361){\line( 1,-1){600}}
 \put( 601, 239){\line( 1, 0){1500}}
 \put( 601,-961){\line( 1, 0){1500}}
 \put(3301, -61){\line( 1, 0){1200}}
 \put(4501, -61){\line( 1, 1){300}}
 \put(3301,-661){\line( 1, 0){1200}}
 \put(4501,-661){\line( 1,-1){300}}
 \put(3601, -61){\line( 1,-1){600}}
 \put(3301, 239){\line( 1, 0){1800}}
 \put(3301,-961){\line( 1, 0){1800}}
 \put( 601, -61){\line( 1,-1){300}}
 \put( 601,-661){\line( 1, 1){300}}
 \put(4501, -61){\line( 1, 1){300}}
 \put(4501,-661){\line( 1,-1){300}}
 \end{picture}
 \]
 This equation can be proved as follows:
 \[
 \setlength{\unitlength}{0.0005in}%
 \begin{picture}(9094,1366)(590,-1044)
  \put(901,-61){\circle*{100}}
  \put(1501,-661){\circle*{100}}
  \put(2101,239){\circle*{100}}
  \put(2101,-961){\circle*{100}}
  \put(3601,-961){\circle*{100}}
  \put(3601,-61){\circle*{100}}
  \put(3901,239){\circle*{100}}
  \put(4501,-961){\circle*{100}}
  \put(2401,239){\circle*{100}}
  \put(2401,-961){\circle*{100}}
  \put(4801,239){\circle*{100}}
  \put(4801,-961){\circle*{100}}
  \put(6001,-61){\circle*{100}}
  \put(6901,-961){\circle*{100}}
  \put(7201,-961){\circle*{100}}
  \put(8401,-361){\circle*{100}}
  \put(8701,-361){\circle*{100}}
  \put(9301,-961){\circle*{100}}
  \put(9301,239){\circle*{100}}
  \put(9601,-961){\circle*{100}}
  \put(9601,239){\circle*{100}}
  \put(6001,-961){\circle*{100}}
  \put(602,-62){\line( 1, 0){1200}}
  \put(1802,-62){\line( 1, 1){300}}
  \put(602,-662){\line( 1, 0){1200}}
  \put(1802,-662){\line( 1,-1){300}}
  \put(902,-62){\line( 1,-1){600}}
  \put(602,238){\line( 1, 0){1800}}
  \put(602,-962){\line( 1, 0){1800}}
  \put(3301,-661){\line( 1,-1){300}}
  \put(3301,-61){\line( 1, 0){300}}
  \put(3601,-61){\line( 1,-1){900}}
  \put(3601,-61){\line( 1, 1){300}}
  \put(3301,239){\line( 1, 0){1500}}
  \put(3301,-961){\line( 1, 0){1500}}
  \put(5701,239){\line( 1,-1){1200}}
  \put(5701,-61){\line( 1, 0){300}}
  \put(5701,-961){\line( 1, 0){1500}}
  \put(8101,-61){\line( 1,-1){300}}
  \put(8401,-361){\line(-1,-1){300}}
  \put(8401,-361){\line( 1, 0){300}}
  \put(9301,239){\line(-1,-1){600}}
  \put(8701,-361){\line( 1,-1){600}}
  \put(8101,239){\line( 1, 0){1500}}
  \put(8101,-961){\line( 1, 0){1500}}
  \put(5701,-661){\line( 1,-1){300}}
  \put(2776,-436){\makebox(0,0){$=$}}
  \put(5176,-436){\makebox(0,0){$=$}}
  \put(7576,-436){\makebox(0,0){$=$}}
 \end{picture}
 \]
 The first equality uses associativity of $\mu$, the second uses the
 defining property of $\psi$, and the third uses the same two ideas
 backwards. 

 We still need to check that $\eta$ and $\psi$ are the \emph{unique}
 maps giving a Frobenius structure.  For $\eta$ this is easy, because
 the unit for a commutative and associative product is always unique.
 For $\psi$, suppose that $\phi\:A\xra{}A\Smash A$ is another map
 giving a Frobenius structure.  We then have the following equations:
 \[
 \setlength{\unitlength}{0.0005in}%
 \begin{picture}(8795,1366)(589,-1044)
 \put(1201,-361){\circle*{100}}
 \put(1801,-961){\circle*{100}}
 \put(1801, 239){\circle*{100}}
 \put(2101, 239){\circle*{100}}
 \put(2101,-961){\circle*{100}}
 \put(3301, -61){\circle*{100}}
 \put(3601, -61){\circle*{100}}
 \put(3901, 239){\circle*{100}}
 \put(4501,-961){\circle*{100}}
 \put(4801,-961){\circle*{100}}
 \put(6001, -61){\circle*{100}}
 \put(6901,-961){\circle*{100}}
 \put(7201,-961){\circle*{100}}
 \put(8701,-361){\circle*{100}}
 \put(8901,-361){\circle*{100}}
 \put(9301,-361){\circle*{100}}
 \put( 601,-361){\line( 1, 0){600}}
 \put(1201,-361){\line( 1, 1){600}}
 \put(1201,-361){\line( 1,-1){600}}
 \put( 601,-961){\line( 1, 0){1500}}
 \put( 601, 239){\line( 1, 0){1500}}
 \put(3001, 239){\line( 1,-1){300}}
 \put(3301, -61){\line(-1,-1){300}}
 \put(3301, -61){\line( 1, 0){300}}
 \put(3601, -61){\line( 1, 1){300}}
 \put(3601, -61){\line( 1,-1){900}}
 \put(3001,-961){\line( 1, 0){1800}}
 \put(5701,-361){\line( 1, 1){300}}
 \put(5701, 239){\line( 1,-1){1200}}
 \put(5701,-961){\line( 1, 0){1500}}
 \put(8101, 239){\line( 1, 0){200}}
 \put(8301, 239){\line( 1,-1){600}}
 \put(8701,-361){\line(-1,-1){600}}
 \put(8101,-361){\line( 1, 0){1200}}
 \put(2506,-436){\makebox(0,0){$=$}}
 \put(5236,-436){\makebox(0,0){$=$}}
 \put(7621,-436){\makebox(0,0){$=$}}
 \put(1066,-616){\makebox(0,0){$\phi$}}
 \put(3556,-391){\makebox(0,0){$\phi$}}
 \end{picture}
 \]
 The first equality is the interchange axiom, the second is the counit
 property of $\ep$, and the third is the associativity of $\mu$.  This
 shows that $\phi$ has the defining property of $\psi$, so $\phi=\psi$
 as required.
\end{proof}
\begin{scholium}\label{sch-maps-duality}
 Let $(A,\eta,\ep,\mu,\psi)$ be a Frobenius object.  Give $A$ the
 inner product $b=\ep\mu$ and give $A\Smash A$ the inner product
 $b\Smash b$.  Then $\eta\:S\xra{}A$ is adjoint to $\ep\:A\xra{}S$ and
 $\psi\:A\xra{}A\Smash A$ is adjoint to $\mu\:A\Smash A\xra{}A$.
\end{scholium}
\begin{proof}
 This is implicit in the proof of the proposition. 
\end{proof}
\begin{scholium}\label{sch-determine-ep}
 The map $\ep\:A\xra{}S$ is the unique one such that
 $(\ep\Smash 1)\psi\eta=\eta\:S\xra{}A$.
\end{scholium}
\begin{proof}
 We saw in the proof of the proposition that $\eta$ is the unique map
 giving the following equality of Penrose diagrams.
 \[
 \setlength{\unitlength}{0.0005in}%
 \begin{picture}(3695,866)(589,-444)
  \put(901,-361){\circle*{100}}
  \put(1501,239){\circle*{100}}
  \put(1501,439){\makebox(0,0){$\mu$}}
  \put(2101,239){\circle*{100}}
  \put(2101,469){\makebox(0,0){$\ep$}}
  \put(4201,239){\circle*{100}}
  \put(4201,469){\makebox(0,0){$\ep$}}
  \put( 601,239){\line( 1, 0){1500}}
  \put(901,-361){\line( 1, 1){600}}
  \put(3301,239){\line( 1, 0){900}}
  \put(2626,-136){\makebox(0,0){$=$}}
  \put(676,-436){\makebox(0,0){$\eta$}}
 \end{picture}
 \]
 The claim follows by working in the opposite category and using
 Remark~\ref{rem-axioms-self-dual}. 
\end{proof}

\begin{remark}\label{rem-A-linear}
 Let $A$ and $B$ be Frobenius objects, and suppose that $f\:A\xra{}B$
 is a ring map with respect to $\eta$ and $\mu$.  We can use this to
 make $B$ into an $A$-module.  We claim that $f^t\:B\xra{}A$ is
 automatically a map of $A$-module objects.  We will give the proof in
 the category of vector spaces over a field; it can easily be made
 diagrammatic.  The claim is that $f^t(f(a)b)=a\,f^t(b)$.  It suffices
 to prove that $b_A(a',f^t(f(a)b))=b_A(a',a\,f^t(b))$.  The left hand
 side is $b_B(f(a'),f(a)b)=\ep_B(f(a)f(a')b)=\ep_B(f(aa')b)$.  The
 right hand side is 
 \[ \ep_A(aa'\,f^t(b))=b_A(aa',f^t(b))=b_B(f(aa'),b)=\ep_B(f(aa')b),\]
 as required.
\end{remark}

\section{The trace form}
\label{sec-trace-form}

We now construct an interesting map $A\xra{}S$ which may or may not be
a Frobenius form.
\begin{definition}\label{defn-trace-form}
 Let $A$ be an arbitrary commutative ring object in an additive
 compact closed category $\CC$.  We can then transpose the
 multiplication map $\mu\:A\Smash A\xra{}A$ to get
 $\rho\:A\xra{}DA\Smash A$ and compose with the evaluation map
 $DA\Smash A=F(A,S)\Smash A\xra{}S$ to get a map $\tht\:A\xra{}S$.
 This is called the \emph{trace form}.  
\end{definition}

\begin{remark}\label{rem-separable}
 If $K$ is a ring and $\CC$ is the category of finitely generated free
 modules over $K$ then $\tht(a)$ is just the trace of the
 multiplication-by-$a$ map.  Now suppose that $K$ is a perfect field.
 One can check that $\tht$ is a Frobenius form if and only if $A$ has
 no nilpotents, if and only if $A$ is a finite product of finite
 extension fields of $K$ (this is well-known and can mostly be
 extracted from~\cite[Section I.1]{frta:ant}, for example).
\end{remark}

\begin{proposition}\label{prop-frobenius-trace}
 Let $A$ be a Frobenius object in an additive compact closed category
 $\CC$.  Then the trace form $\tht$ is given by
 $\tht=b\psi=\ep\mu\psi\:A\xra{}S$.  Moreover, if we define
 $\al:=\mu c=\mu\psi\eta\:S\xra{}A$ then $\tht=b(\al\Smash 1)$.
\end{proposition}
\begin{proof}
 The adjunction between the functors $A\Smash(-)$ and
 $F(A,-)=DA\Smash(-)$ is given by two maps
 $\text{unit}\:S\xra{}DA\Smash A$ and
 $\text{eval}\:DA\Smash A\xra{}S$.  It follows from the basic theory
 of pairings and duality~\cite[Chapter III]{lemast:esh} that the
 following diagrams commute:
 \begin{diag}
  \node{A\Smash A} 
  \arrow{s,l}{b^\#\Smash 1}
  \arrow{se,t}{b}
  \node[2]{S}
  \arrow{se,b}{\text{unit}}
  \arrow{e,t}{c}
  \node{A\Smash A}
  \arrow{s,r}{b^\#\Smash 1} \\
  \node{DA\Smash A}
  \arrow{e,b}{\text{eval}}
  \node{S}
  \node[2]{DA\Smash A}
 \end{diag}
 It follows that the following diagram commutes:
 \begin{diag}
  \node{A}
  \arrow{s,=}
  \arrow{e,t}{c\Smash 1}
  \node{A\Smash A\Smash A}
  \arrow{s,l}{b^\#\Smash 1\Smash 1}
  \arrow{e,t}{1\Smash\mu}
  \node{A\Smash A}
  \arrow{s,r}{b^\#\Smash 1}
  \arrow{e,t}{b}
  \node{S} 
  \arrow{s,=} \\
  \node{A}
  \arrow{e,b}{\text{unit}\Smash 1}
  \node{DA\Smash A\Smash A}
  \arrow{e,b}{1\Smash \mu}
  \node{DA\Smash A}
  \arrow{e,b}{\text{eval}}
  \node{S}
 \end{diag}
 On the bottom row, the composite of the first two maps is $\rho$ so
 the whole composite is just $\tht$.  Thus,
 $\tht=b(1\Smash\mu)(c\Smash 1)$.  To complete the proof, it is
 easiest to think in terms of TQFT's.  Let $M$ be a torus with a small
 open disc removed.  We leave it to the reader to check that
 $\al=\mu c$ is represented by $M$, considered as a cobordism from
 $\emptyset$ to $S^1$.  Moreover, the maps $b(1\Smash\mu)(c\Smash 1)$,
 $b\psi$ and $b(\alpha\Smash 1)$ are all represented by $M$ considered
 as a cobordism from $S^1$ to $\emptyset$.  The proposition follows.
\end{proof}

\section{Manifolds}
\label{sec-manifolds}

We next show how to use manifolds to construct Frobenius objects in
suitable categories of module spectra.  This is of course just a
reformulation of Atiyah-Poincar\'{e} duality, but it is a nice
illustration of the theory of Frobenius objects.  It is also
strikingly formally similar to the constructions in the $K(n)$-local
stable category which we discuss later.

Let $\CM$ be the category of even-dimensional closed manifolds $M$
equipped with a complex structure on the stable normal bundle, or
equivalently a complex orientation on the map from $M$ to the
one-point manifold; we refer to Quillen's work~\cite{qu:epc} for a
careful discussion of what this means.  

Next, let $MP$ denote the Thom spectrum of the tautological virtual
complex bundle over $\Z\tm BU$, so that
$MP=\bigWedge_{n\in\Z}\Sg^{2n}MU$ and $\Sg^2MP\simeq MP$.  More
generally, if $V$ is a complex bundle over a space $X$ then there is a
canonical Thom class $u_V\:X^V\xra{}MP$ which combines with the usual
diagonal map $X^V\xra{}X_+\Smash X^V$ to give a canonical equivalence
$MP\Smash X^V\simeq MP\Smash\Sgi X_+$.  With a little care, this also
goes through for virtual bundles. 

The spectrum $MP$ can be constructed as an $E_\infty$ ring spectrum,
and thus as a strictly commutative ring spectrum (or ``$S$-algebra'')
in the EKMM category~\cite{ekmm:rma}.  We can thus define a category
of $MP$-modules in the strict sense, and the associated derived
category $\DD=\DD_{MP}$.  (There are also other approaches to our
results using less technology.)  The category $\DD$ is a unital
algebraic stable homotopy category in the sense of~\cite{hopast:ash};
in particular it is a closed symmetric monoidal category.  We write
$\CF$ for the thick subcategory of $\DD$ generated by $MP$, which is
the same as the category of small or strongly dualisable
objects~\cite[Theorem 2.1.3(d)]{hopast:ash}.  This is clearly a
compact closed category.

Define $T\:\CM\xra{}\CF$ by $T(M)=MP\Smash\Sgi M_+$.  This is clearly
a covariant functor that converts products to smash products and
disjoint unions to wedges.

Now suppose we have a smooth map $f\:M\xra{}N$ of closed manifolds.
Let $j\:M\xra{}\R^k$ be a smooth map such that 
$(j,f)\:M\xra{}\R^k\tm N$ is a closed embedding, with normal bundle
$\nu_{(j,f)}$ say.  This is stably equivalent to $k+\nu_M-f^*\nu_N$.
The Pontrjagin-Thom construction applied to the embedding $(j,f)$
gives a map $\Sg^kN_+\xra{}M^{\nu_{(j,f)}}$ and thus a stable map
$f^!\:\Sgi N_+\xra{}M^{\nu_M-f^*\nu_N}$.

Now suppose that $M$ and $N$ have specified complex orientations, so
they are objects of $\CM$.  Then the virtual bundle
$\nu_f=\nu_M-f^*\nu_N$ has a canonical complex structure, so there is
a canonical equivalence $MP\Smash M^{\nu_M-f^*\nu_N}\simeq T(M)$.
Thus, by smashing $f^!$ with $MP$ we get a map $Uf\:T(N)\xra{}T(M)$.
One can check that this construction gives a contravariant functor
$U\:\CM\xra{}\DD$, which again converts products to smash products and
disjoint unions to wedges.  If $f$ is a diffeomorphism, one checks
easily that $U(f)=T(f)^{-1}$.  We also have the following ``Mackey
property''.  Suppose we have a commutative square in $\CM$:
\begin{diag}
 \node{K} \arrow{e,t}{f} \arrow{s,l}{g} \node{L} \arrow{s,r}{h} \\
 \node{M} \arrow{e,b}{k}                \node{N.}
\end{diag}
Suppose also that the square is a pullback and the maps $h$ and $k$
are transverse to each other, so that when $x\in K$ with
$hf(x)=kg(x)=y$ say, the map of tangent spaces 
\[ (Dh,Dk)\:T_{f(x)}L \oplus T_{g(x)}M \xra{} T_yN \] 
is surjective.  We then have $U(h)T(k)=T(f)U(g)$, as one sees directly
from the geometry.

For any manifold $M\in\CM$, there is a unique map $\ep\:M\xra{}1$,
where $1$ is the one-point manifold.  We also have a diagonal map
$\psi\:M\xra{}M\tm M$.  We allow ourselves to write $\ep$ and $\psi$
for $T(\ep)$ and $T(\psi)$, and we also write $\eta=U(\ep)$ and
$\mu=U(\psi)$.  

\begin{proposition}\label{prop-manifolds-frobenius}
 The above maps make $T(M)$ into a Frobenius object in $\CM$.  If we
 use the resulting inner product, then for any map $f\:M\xra{}N$ in
 $\CM$ we have $T(f)^t=U(f)$.
\end{proposition}
\begin{proof}
 If we make $\CM$ into a symmetric monoidal category using the
 cartesian product, it is clear that $\ep$ and $\psi$ make $M$ into a
 comonoid object.  As $T$ and $U$ are monoidal functors, the first
 covariant and the second contravariant, we see that
 $\Sgi M_+=T(M)=U(M)$ is a monoid object under $\mu$ and $\eta$, and a
 comonoid object under $\psi$ and $\ep$.  For the interchange axiom
 $\psi\mu=(1\Smash\mu)(\psi\Smash 1)$, we note that the following
 diagram is a transverse pullback and apply the Mackey property.
 \begin{diag}
  \node{M}           \arrow{e,t}{\psi}      \arrow{s,l}{\psi}
  \node{M\tm M}                            \arrow{s,r}{1\tm\psi}\\
  \node{M\tm M}      \arrow{e,b}{\psi\tm 1} 
  \node{M\tm M\tm M.}
 \end{diag}
 Similarly, to prove that $T(f)^t=U(f)$, we note that the following
 square is a transverse pullback:
 \begin{diag}
  \node{M}      \arrow{e,t}{f}      \arrow{s,l}{(1,f)}
  \node{N}                          \arrow{s,r}{\psi}    \\
  \node{M\tm N} \arrow{e,b}{f\tm 1} \node{N\tm N.}
 \end{diag}
 We then apply the Mackey property, noting that $(1,f)=(1\tm f)\psi_M$;
 this gives the following commutative diagram:
 \begin{diag}
  \node{T(M)\Smash T(N)}
  \arrow[2]{e,t}{T(f)\Smash 1}
  \arrow{s,l}{1\Smash U(f)}
  \node[2]{T(N)\Smash T(N)}
  \arrow{s,r}{\mu_N}            \\
  \node{T(M)\Smash T(M)}
  \arrow{e,b}{\mu_M}
  \node{T(M)}
  \arrow{e,b}{T(f)}
  \node{T(N).}
 \end{diag}
 We then compose with $\ep_N$, noting that $\ep_NT(f)=\ep_M$ and
 $\ep_M\mu_M=b_M$ and $\ep_N\mu_N=b_N$.  We conclude that
 $b_N(T(f)\Smash 1)=b_M(1\Smash U(f))$, so $T(f)^t=U(f)$ as claimed.
\end{proof}

\section{Groupoids}
\label{sec-groupoids}

Let $\CG$ denote the category of groupoids and functors between them,
and let $\CGb$ be the quotient category in which two functors are
identified if there is a natural isomorphism between them.  We say
that a groupoid $G$ is \emph{finite} if there are only finitely many
isomorphism classes of objects, and $G(a,b)$ is finite for any 
$a,b\in G$.  We write $\CGf$ for the category of finite groupoids, and
$\CGbf$ for the obvious quotient category.

We next exhibit an equvalence between $\CGb$ and a certain homotopy
category of spaces.  As usual in homotopy theory, it will be
convenient to work with compactly generated weakly Hausdorff spaces
(so we have Cartesian closure).  Let $\CB$ be the category of such
spaces $X$ for which $\pi_k(X,x)=0$ for all $k>1$ and all $x\in X$.
We also write $\CBb$ for the associated homotopy category (in which
weak equivalences are inverted), and we let $\CBf$ and $\CBbf$ be the
subcategories whose objects are those $X\in\CB$ for which $\pi_0X$ is
finite and $\pi_1(X,x)$ is finite for any basepoint $x$.

Milgram's classifying space construction gives a functor
$B\:\CG\xra{}\CB$.  One can also define a functor
$\Pi_1\:\CB\xra{}\CG$: the set of objects of $\Pi_1(X)$ is $X$, and
the set of morphisms from $x$ to $y$ is the set of paths from $x$ to
$y$ modulo homotopy relative to the endpoints.  Both $\CG$ and $\CB$
have finite products and coproducts, and both our functors preserve
them.

It is easy to check that these constructions give equivalences
$\CGb\simeq\CBb$ and $\CGbf\simeq\CBbf$.

Any (finite) group $G$ can be regarded as a (finite) groupoid with one
object.  If $G$ and $H$ are groups then $\CG(G,H)$ is the set of
homomorphisms from $G$ to $H$, and $\CGb(G,H)$ is the set of conjugacy
classes of such homomorphisms.

Conversely, if $G$ is a finite groupoid then we can choose a family
$\{a_i\}_{i\in I}$ containing precisely one object of $G$ from each
isomorphism class and then let $H_i$ be the group $G(a_i,a_i)$.  We
find that $G\simeq\coprod_iH_i$ in $\CGb$.  Thus, all our questions
about groupoids can be reduced to questions about groups by some
unnatural choices.  Our next lemma sharpens this slightly.

\begin{definition}\label{defn-indiscrete}
 A groupoid $G$ is \emph{discrete} if all its maps are identity maps,
 and \emph{indiscrete} if there is precisely one map from $a$ to $a'$
 for all $a,a'\in G$.  
\end{definition}
\begin{remark}\label{rem-indiscrete}
 The category of discrete groupoids is equivalent to that of sets, as
 is the category of indiscrete groupoids.  The classifying space of a
 discrete groupoid is discrete, and that of a nonempty indiscrete
 groupoid is contractible.
\end{remark}

\begin{lemma}\label{lem-connected-split}
 Any nonempty connected groupoid is isomorphic to $A\tm H$ for some
 nonempty indiscrete groupoid $A$ and some group $H$.  Thus, any
 groupoid is isomorphic to $\coprod_IA_i\tm H_i$ for some family of
 nonempty indiscrete groupoids $A_i$ and groups $H_i$.
\end{lemma}
\begin{proof}
 Let $G$ be a connected groupoid.  Choose an object $x\in G$ and let
 $H$ be the group $G(x,x)$.  Let $A$ be the indiscrete groupoid with
 $\obj(A)=\obj(G)$, and for each $a\in A$ choose a map $k_a\:x\xra{}a$
 in $G$.  Put $B=A\tm H$, so $\obj(B)=\obj(G)$ and $B(a,a')=H$ for all
 $a,a'$.  Composition is given by multiplication in $H$.  Define
 $u\:B\xra{}G$ by $u(a)=a$ on objects, and 
 \[ u_{a,a'}(h) = (a \xra{k_a^{-1}} x \xra{h} x \xra{k_{a'}} a') \]
 on morphisms.  This is easily seen to be functorial and to be an
 isomorphism.  

 The generalisation to the disconnected case is immediate.
\end{proof}

\subsection{Model category structure}

We now complete an exercise assigned by Anderson~\cite{an:fgr} to his
readers, by verifying that his definitions (reproduced below) do
indeed make the category $\CG$ into a closed model category in the
sense of Quillen~\cite{qu:ha} (see also~\cite{dwsp:htm} for an
exposition and survey of more recent literature).  As well as being
useful for our applications, this seems pedagogically valuable, as the
verification of the axioms is simpler than in most other examples.
The homotopy theory of the category of all small categories has been
extensively studied (see~\cite{qu:haki} for example), but the case of
groupoids is easier so it makes sense to treat it independently.

\begin{definition}
 We say that a functor $u\:G\xra{}H$ of groupoids is
 \begin{itemize}
  \item[(a)] a \emph{weak equivalence} if it is full, faithful and
   essentially surjective (in other words, an equivalence of
   categories);
  \item[(b)] a \emph{cofibration} if it is injective on objects; and
  \item[(c)] a \emph{fibration} if for all $a\in G$, $b\in H$ and
   $h\:u(a)\xra{}b$ there exists $g\:a\xra{}a'$ in $G$ such that
   $u(a')=b$ and $u(g)=h$.
 \end{itemize}
 As usual, an \emph{acyclic fibration} means a fibration that is also
 an equivalence, and similarly for acyclic cofibrations.
\end{definition}

\begin{remark}\label{rem-hom-fibration}
 Let $u\:G\xra{}H$ be a homomorphism of groups.  Then $u$ is
 automatically a cofibration of groupoids, and it is a fibration iff
 it is surjective.  It is an equivalence of groupoids iff it is an
 isomorphism.
\end{remark}
\begin{remark}
 Let $v\:X\xra{}Y$ be a map of sets.  If we regard $X$ and $Y$ as
 discrete categories then $v$ is automatically a fibration.  It is a
 cofibration iff it is injective, and an equivalence iff it is
 bijective.  If we regard $X$ and $Y$ as indiscrete categories then
 $v$ is automatically an equivalence (unless $\emptyset=X\neq Y$).  It
 is a cofibration iff it is injective, and a fibration iff it is
 surjective. 
\end{remark}

\begin{theorem}
 The above definitions make $\CG$ into a closed model category.
\end{theorem}
\begin{proof}
 We need to verify the following axioms, numbered as
 in~\cite{dwsp:htm}:
 \begin{itemize}
  \item[MC1:] $\CG$ has finite limits and colimits.
  \item[MC2:] If we have functors $G\xra{u}H\xra{v}K$ and two of $u$,
   $v$ and $vu$ are weak equivalences then so is the third.
  \item[MC3:] Every retract of a weak equivalence is a weak
   equivalence, and similarly for fibrations and cofibrations.
  \item[MC4:] Cofibrations have the left lifting property for acyclic
   fibrations, and acyclic cofibrations have the left lifting property
   for all fibrations.
  \item[MC5:] Any functor $u$ has factorisations $u=pi=qj$ where $i$
   and $j$ are cofibrations, $p$ and $q$ are fibrations, and $i$ and
   $q$ are equivalences.
 \end{itemize}

 MC1: This follows from the fact that $\CG$ is the category of models
 for a left-exact sketch~\cite[Section 4.4]{bawe:ttt}.  More
 concretely, for limits we just have
 $\obj(\invlim_iG_i)=\invlim_i\obj(G_i)$ and
 $\mor(\invlim_iG_i)=\invlim_i\mor(G_i)$.  Similarly, for coproducts
 we have $\obj(\coprod_iG_i)=\coprod_i\obj(G_i)$ and
 $\mor(\coprod_iG_i)=\coprod_i\mor(G_i)$.  Coequalisers are more
 complicated and best handled by the adjoint functor theorem.

 MC2: This is easy.

 MC3: Let $v$ be an equivalence and let $u\:G\xra{}H$ be a retract of
 $v$.  Then $\pi_0(u)$ is a retract of $\pi_0(v)$, so $\pi_0(u)$ is a
 bijection and so $u$ is essentially surjective.  If $a,b\in G$ then
 $u_{a,b}\:G(a,b)\xra{}H(ua,ub)$ is a retract of a map of the form
 $v_{c,d}$ and thus is a bijection, so $u$ is full and faithful.  Thus
 $u$ is an equivalence as required.

 It is clear that a retract of a cofibration is a cofibration.  

 For fibrations, let $1$ be the terminal groupoid.  Let $I$ be the
 groupoid with objects $\{0,1\}$ and two non-identity morphisms
 $u\:0\xra{}1$ and $u^{-1}\:1\xra{}0$.  Let $i\:1\xra{}I$ be the
 inclusion of $\{0\}$.  Then fibrations are precisely the maps with
 the right lifting property for $i$, and it follows that a retract of
 a fibration is a fibration.

 MC4: Consider a commutative square as follows, in which $i$ is a
 cofibration and $p$ is a fibration.
 \begin{diag}
  \node{G} \arrow{s,l}{i} \arrow{e,t}{u} \node{K} \arrow{s,r}{p} \\ 
  \node{H}                \arrow{e,b}{v} \node{L.}
 \end{diag}
 Because $p$ is a fibration, it is easy to see that the image of $p$
 is replete: if $d\in L$ is isomorphic to $pc$ then $d$ has the form
 $pc'$ for some $c'\in K$.

 Suppose that $p$ is an equivalence; we must construct a functor
 $w\:H\xra{}K$ such that $pw=v$ and $wi=u$.  As $p$ is essentially
 surjective and the image is replete, we see that $\obj(p)$ is
 surjective.  By assumption $i$ is a cofibration so $\obj(i)$ is
 injective.  Define a map $w\:\obj(H)\xra{}\obj(K)$ by putting
 $w(i(a))=u(a)$ for $a\in\obj(A)$ and choosing $w(b)$ to be any
 preimage under $p$ of $v(b)$ if $b\not\in\image(i)$.  Clearly $pw=v$
 and $wi=u$ on objects.  Given $b,b'\in H$ we define $w_{b,b'}$ to be
 the composite
 \[ H(b,b') \xra{v_{b,b'}} L(vb,vb') = L(pwb,pwb')
     \xra{p_{wb,wb'}^{-1}} K(wb,wb').
 \]
 One can check that this makes $w$ a functor with $pw=v$.  Also
 $pwi=vi=pu$ on morphisms and $wi=u$ on objects and $p$ is faithful;
 it follows that $wi=u$ on morphisms, as required.

 Now remove the assumption that $p$ is an equivalence, and suppose
 instead that $i$ is an equivalence.  We must again define a functor
 $w\:H\xra{}K$ making everything commute.  As $i$ is injective on
 objects we can choose $r\:\obj(H)\xra{}\obj(G)$ with $ri=1$.  As
 $\pi_0(i)$ is a bijection we find that $\pi_0(r)=\pi_0(i)^{-1}$ so we
 can choose isomorphisms $\eta_b\:b\xra{}ir(b)$ for all $b\in H$.  If
 $b=i(a)$ for some (necessarily unique) object $a$, we have $rb=a$ and
 $irb=b$, and we choose $\eta_b=1_b$ in this case.  There is a unique
 way to make $r$ a functor $H\xra{}G$ such that $\eta$ is natural:
 explicitly, the map $r_{b,b'}$ is the composite
 \[ H(b,b') \xra{\eta_*\eta^*} H(irb,irb') 
    \xra{i_{rb,rb'}^{-1}} G(rb,rb').
 \]
 Next, if $b\in\image(i)$ we define $wb=urb$ and
 $\zt_b=1\:wb\xra{}urb$.  If $b\not\in\image(i)$ we instead apply the
 fibration axiom for $p$ to the map $v\eta_b\:vb\xra{}virb=purb$ to
 get an object $wb\in K$ and a morphism $\zt_b\:wb\xra{}urb$ such that
 $pwb=vb$ and $p\zt_b=v\eta_b$.  Note that these last two equations
 also hold in the case $b\in\image(i)$.  There is a unique way to make
 $w$ into a functor such that $\zt\:w\xra{}ur$ is natural.  Clearly
 $wi=u$ as functors, and $pw=v$ on objects.  Given $h\:b\xra{}b'$ in
 $H$ we can apply $p$ to the naturality square for $\zt$ and then use
 the naturality of $\eta$ to deduce that $pwh=vh$; thus $pw=v$ on
 morphisms, as required.

 MC5: Consider a functor $u\:G\xra{}H$.  Let $K$ be the category whose
 objects are triples $(a,b,k)$, with $a\in G$ and $b\in H$ and
 $k\:u(a)\xra{}b$.  The morphisms from $(a,b,k)$ to $(a',b',k')$ are
 the pairs $(g,h)$ where $g\:a\xra{}a'$ and $h\:b\xra{}b'$ and the
 following diagram commutes:
 \begin{diag}
  \node{u(a)}  \arrow{s,l}{k} \arrow{e,t}{u(g)}
  \node{u(a')} \arrow{s,r}{k'}                  \\
  \node{b}     \arrow{e,b}{h} \node{b'.}
 \end{diag}
 We also consider the category $L$ with the same objects as $K$, but
 with $L(a,b,k;a',b',k')=H(b,b')$, so there is an evident functor
 $v\:K\xra{}L$.  There is also a functor $i\:G\xra{}K$ given by
 $i(a)=(a,ua,1_{ua})$ and a functor $q\:L\xra{}H$ given by
 $q(a,b,k)=b$; we put $j=vi$ and $p=qv$.  It is clear that
 $u=qvi=qj=pi$ and that $i$ and $j$ are cofibrations and that $i$ is
 full and faithful.  If $(a,b,k)\in K$ then
 $(1_a,k)\:i(a)\xra{}(a,b,k)$ so $i$ is essentially surjective and
 thus an equivalence.  The functor $q$ is clearly full and faithful,
 and its image is the repletion of the image of $u$.  

 We next claim that $p$ and $q$ are fibrations.  Suppose that
 $(a,b,k)\in\obj(K)$ and $h\:b=q(a,b,k)\xra{}b'$.  Then
 $(a,b',hk)\in\obj(K)$ and $(1_a,h)\:(a,b,k)\xra{}(a,b',hk)$ and
 $q(1_a,h)=h$.  This shows that $q$ is a fibration, and the same
 construction also shows that $p$ is a fibration.  

 We now have a factorisation $u=pi$ as required by axiom MC3.  
 If $u$ is essentially surjective then the same is true of $q$ and
 thus $q$ is an equivalence and so the factorisation $u=qj$ is also as
 required.  If $u$ is not essentially surjective then we let $L'$ be the
 full subcategory of $H$ consisting of objects not in the repletion of
 the image of $u$ and let $q'\:L'\xra{}H$ be the inclusion.  We then
 have an acyclic fibration $(q,q')\:L\amalg L'\xra{}H$ and a
 cofibration $G\xra{j}L\xra{}L\amalg L'$ whose composite is $u$, as
 required. 
\end{proof}

\begin{proposition}\label{prop-coproper}
 The above model category structure is right proper (in other words,
 the pullback of a weak equivalence along a fibration is a weak
 equivalence.)
\end{proposition}

\begin{proof}
 Consider a pullback square as follows, in which $v$ is a weak
 equivalence and $q$ is a fibration. 
 \begin{diag}
  \node{G} \arrow{s,l}{p} \arrow{e,t}{u} \node{K} \arrow{s,r}{q} \\ 
  \node{H}                \arrow{e,b}{v} \node{L.}
 \end{diag}
 Suppose that $a,a'\in G$ and put
 $d=qu(a)=vp(a)$ and $d'=qu(a')=vp(a')$.  By the construction of
 pullbacks in $\CG$, we see that the following square is a pullback
 square of sets:
 \begin{diag}
  \node{G(a,a')}
  \arrow{s,l}{p_{a,a'}}
  \arrow{e,t}{u_{a,a'}}
  \node{K(u(a),u(a'))}
  \arrow{s,r}{q_{u(a),u(a')}} \\ 
  \node{H(p(a),p(a'))}
  \arrow{e,b}{v_{p(a),p(a')}}
  \node{L(d,d').}
 \end{diag}
 As $v$ is a weak equivalence, the map $v_{p(a),p(a')}$ is a
 bijection, and it follows that the same is true of $u_{a,a'}$.  This
 means that $u$ is full and faithful.  

 Next suppose we have $c\in K$, so $q(c)\in L$.  As $v$ is essentially
 surjective there exists $b\in H$ and $l\:q(c)\xra{}v(b)$ in $L$.  As
 $q$ is a fibration there is a map $k\:c\xra{}c'$ in $K$ with
 $q(c')=v(b)$ and $q(k)=l$.  By the pullback property there is a
 unique $a\in G$ with $u(a)=c'$ and $p(a)=b$.  Thus $u(a)\simeq c$,
 proving that $u$ is essentially surjective and thus an equivalence.
\end{proof}

\subsection{Classifying spaces}

Let $N$ be the nerve functor from groupoids to simplicial sets, 
and put $BG=|NG|$; this is called the \emph{classifying space} of
$G$.  It is easy to see that $N$ converts groupoids to Kan complexes
and fibrations to Kan fibrations, and that it preserves coproducts and
finite limits.  The geometric realisation functor preserves coproducts
(easy) and finite limits~\cite[Theorem 4.3.16]{frpi:cst} and it
converts Kan fibrations to fibrations~\cite{qu:grk} (see
also~\cite[Theorem 4.5.25]{frpi:cst}).  Thus, the composite functor
$B\:\CG\xra{}\CB$ preserves coproducts, finite limits and fibrations.

\subsection{Homotopy pullbacks}

\begin{definition}
 Suppose we have functors $G\xra{u}H\xla{v}K$ of groupoids.  We define
 a new groupoid $L$ whose objects are triples $(a,c,h)$ with $a\in G$
 and $c\in K$ and $h\:u(a)\xra{}v(c)$.  The morphisms from $(a,c,h)$
 to $(a',c',h')$ are pairs $(r,s)$ where $r\:a\xra{}a'$ and
 $s\:c\xra{}c'$ and the following diagram commutes:
 \begin{diag}
  \node{u(a)} \arrow{e,t}{u(r)} \arrow{s,l}{h}
  \node{u(a')}                  \arrow{s,r}{h'} \\
  \node{v(c)} \arrow{e,b}{v(s)} \node{v(c').}
 \end{diag}
 We also define functors $K\xla{u'}L\xra{v'}G$ by $u'(a,c,h)=a$ and
 $v'(a,c,h)=c$, and a natural transformation $\phi\:uv'\xra{}vu'$ by
 $\phi_{(a,c,h)}=h$.  This gives a square as follows, which commutes
 in $\CGb$:
 \begin{diag}
  \node{L} \arrow{e,t}{v'} \arrow{s,l}{u'} 
  \node{G}                 \arrow{s,r}{u} \\
  \node{K} \arrow{e,b}{v}  \node{H.}
 \end{diag}
 We call $L$ the \emph{homotopy pullback} of $u$ and $v$.  We say that
 an arbitrary commutative square in $\CGb$ is
 \emph{homotopy-cartesian} if it is isomorphic to one of the above
 form.
\end{definition}

\begin{remark}\label{rem-actual-pullback}
 We can also consider the actual pullback rather than the homotopy
 pullback, which can be identified with the full subcategory $M\sse L$
 consisting of pairs $(a,c,1)$ where $u(a)=v(c)$.  One checks that the
 inclusion $M\xra{}L$ is an equivalence if $u$ or $v$ is a fibration.
\end{remark}
\begin{remark}\label{rem-groups-pullback}
 Suppose that $H$ is a group and $u$ and $v$ are inclusions of
 subgroups.  Then $M$ is the group $G\cap K$.  Let $T\sse H$ be a set
 containing one element of each double coset in $G\setminus H/K$; we
 may as well assume that $1\in T$.  We find that $L$ is equivalent to
 the groupoid $\coprod_T G^t\cap K$, and the term indexed by $t=1$ is
 just $M$.  It follows that the map $M\xra{}L$ is an equivalence if
 and only if $H=GK$.  Note that this is only predicted by the previous
 remark when $G=H$ or $K=H$.
\end{remark}
\begin{remark}\label{rem-replace-square}
 By standard methods of abstract homotopy theory, we see that a
 square $S$ in $\CGb$ is homotopy-cartesian iff there is a pullback
 square $S'$ in $\CG$ whose maps are fibrations, which becomes
 isomorphic to $S$ in $\CGb$.  
\end{remark}
\begin{remark}\label{rem-finite-pullback}
 It is easy to see that if $G$, $H$ and $K$ are finite then so is
 their homotopy pullback.
\end{remark}

\begin{definition}
 Suppose we have functors $u,v,s,t$ such that the following square is
 commutative in $\CGb$.
 \begin{diag}
  \node{F} \arrow{e,t}{t} \arrow{s,l}{s} \node{G} \arrow{s,r}{u} \\
  \node{K} \arrow{e,b}{v}                \node{H}
 \end{diag}
 Let $L$ be the homotopy pullback of $u$ and $v$, and let $u',v'$ be
 as above.  Choose an isomorphism $\sg\:ut\xra{}vs$.  We can then
 define a functor $\hat{\sg}\:F\xra{}L$ by
 $\hat{\sg}(d)=(t(d),s(d),\sg_d)$; this has $u'\hat{\sg}=t$ and
 $v'\hat{\sg}=s$.  If $\zt\:s\xra{}s'$ and $\xi\:t'\xra{}t$ and
 $\sg'=v(\zt)\circ s\circ u(\xi)$ then it is easy to see that
 $\sg\simeq\widehat{\sg'}$. 
\end{definition}

\begin{lemma}\label{lem-homotopy-cartesian}
 A square as in the above definition is homotopy Cartesian if and only
 if there exists $\sg\:ut\xra{}vs$ such that $\hat{\sg}\:F\xra{}L$ is
 an equivalence.  
\end{lemma}
\begin{proof}
 If there exists such a map $\sg$ then the square is visibly
 equivalent in $\CGb$ to a homotopy pullback square, and thus is
 homotopy cartesian.  For the converse, suppose that the square is
 homotopy Cartesian.  We can then find a diagram as follows which
 commutes in $\CGb$, such that the outer square is a homotopy
 pullback, and the diagonal functors are equivalences. 
 \begin{diag}
  \node{L_1}
  \arrow[3]{e,t}{v'_1}
  \arrow[3]{s,l}{u'_1}
  \node[3]{G_1} 
  \arrow[3]{s,r}{u_1} \\
  \node[2]{F}
  \arrow{e,t}{t}
  \arrow{s,l}{s}
  \arrow{nw,t}{\dl}
  \node{G}
  \arrow{ne,t}{\al}
  \arrow{s,r}{u} \\
  \node[2]{K}
  \arrow{sw,b}{\gm}
  \arrow{e,b}{v}
  \node{H}
  \arrow{se,b}{\bt} \\
  \node{K_1}
  \arrow[3]{e,b}{v_1}
  \node[3]{H_1}
 \end{diag}
 There is a ``tautological'' natural isomorphism
 $\phi_1\:u_1v_1'\xra{}v_1u_1'$, and we write
 $\rho=\phi_1\dl\:u_1v'_1\dl\xra{}v_1u'_1\dl$ so that
 $\dl=\hat{\rho}$.  As the top and left-hand regions of the diagram
 commute in $\CGb$, we have natural maps $\al t\xra{}v'_1\dl$
 and $u_1'\dl\xra{}\gm s$, which we can use to form a natural map
 \[ \kp = (u_1\al t\xra{}u_1v'_1\dl\xra{\rho}
           v_1u'_1\dl\xra{}v_1\gm s).
 \]
 Using the remark in the preceeding definition, we see that
 $\hat{\kp}\simeq\hat{\rho}=\dl\:F\xra{}L_1$.  As $\dl$ is an
 equivalence, we see that the same is true of $\hat{\kp}$.  Next, we
 note that the functors $u_1\al v',v_1\gm u'\:L\xra{}L'$ are joined by
 the natural map
 \[ \tau=(u_1\al v'\xra{}\bt uv'\xra{\bt(\phi)}
          \bt v u' \xra{} v_1\gm u'),
 \]
 where the first and third maps come from the commutativity of the
 right-hand and bottom regions of the diagram.  This gives a functor
 $\hat{\tau}\:L\xra{}L'$; we leave it to the reader to check directly
 that this is an equivalence.

 Next, consider the composite
 \[ \bt u t \xra{} u_1\al t \xra{\kp} v_1\gm s\xra{}\bt v s. \]
 As $\bt$ is full and faithful, this composite has the form $\bt(\sg)$
 for a unique natural map $\sg\:ut\xra{}vs$, which gives rise to
 $\hat{\sg}\:F\xra{}L$.  One checks directly that
 $\hat{\tau}\hat{\sg}=\hat{\kp}$, and both $\hat{\tau}$ and
 $\hat{\kp}$ are equivalences, so $\hat{\sg}$ is an equivalence, as
 required. 
\end{proof}

\subsection{Coverings and quasi-coverings}

\begin{definition}\label{defn-quasi-covering}
 A functor $u\:G\xra{}H$ is a \emph{covering} if for each $a\in G$ and
 each $h\:u(a)\xra{}b$ in $H$ there is a \emph{unique} pair $(a',g)$
 with $a'\in G$ and $g\:a\xra{}a'$ such that $u(a')=b$ and $u(g)=h$.
 More generally, we say that $u$ is a \emph{quasi-covering} if it can
 be factored as an equivalence followed by a covering.  
\end{definition}
\begin{remark}\label{rem-coverings-preserved}
 It is easy to check that pullbacks, products and composites of
 coverings are coverings.
\end{remark}
\begin{remark}\label{rem-group-coverings}
 A group homomorphism is only a covering if it is an isomorphism.  We
 will see later that it is a quasi-covering iff it is injective.
\end{remark}

\begin{definition}\label{defn-reflects-id}
 A functor $u\:G\xra{}H$ \emph{reflects identities} if whenever
 $g\:a\xra{}a'$ and $u(g)=1_b$ for some $b$, we have $a=a'$ and
 $g=1_a$.  Such a functor is easily seen to be faithful.
\end{definition}

We leave the following easy lemma to the reader.
\begin{lemma}\label{lem-reflects-id}
 A functor $u\:G\xra{}H$ is a covering iff it reflects identities and
 is a fibration.  \qed
\end{lemma}

\begin{proposition}\label{prop-B-covering}
 If $u\:G\xra{}H$ is a covering, then $Bu\:BG\xra{}BH$ is a covering
 map of topological spaces.
\end{proposition}
\begin{proof}
 Suppose for the moment that $H$ is indiscrete and $G$ is connected.
 Then for $a,a'\in G$ we have $G(a,a')\neq\emptyset$ and
 $u\:G(a,a')\xra{}H(ua,ua')$ is injective but the codomain has only
 one element so the same is true of $G(a,a')$.  Thus $u$ is full and
 faithful.  It is also a fibration and $H$ is connected so it is
 surjective on objects.  If $ua=ua'$ then the unique map $a\xra{}a'$
 in $G$ must become an identity map in $H$ but $u$ reflects identities
 so $a=a'$.  We now see that $u$ is an isomorphism so $Bu$ is a
 homeomorphism and thus certainly a covering.

 If $H$ is indiscrete and $G$ is disconnected, we can still show that
 $Bu$ is a covering by looking at one component at a time.

 Now suppose merely that $H$ is connected.  We can then split $H$ as
 $A\tm K$, where $A$ is indiscrete and $K$ is a group, as in
 Lemma~\ref{lem-connected-split}.  Let $K'$ be the indiscrete category
 with object set $K$, and define $q\:K'\xra{}K$ by sending the unique
 morphism $k\xra{}k'$ to $k'k^{-1}\in\mor(K)$.  One checks that
 $BK'=EK$ and that $Bq\:EK\xra{}BK=EK/K$ is the usual covering map.
 Thus, $H'=A\tm K'$ is indiscrete and $r=1\tm q\:H'\xra{}H$ is a
 covering with the property that $Br$ is also a covering.  Now form a
 pullback square as follows:
 \begin{diag}
  \node{G'} \arrow{s,l}{r'} \arrow{e,t}{u'} \node{H'} \arrow{s,r}{r}\\
  \node{G}                  \arrow{e,b}{u}  \node{H.}
 \end{diag}
 Note that $u'$ is a covering.  As $H'$ is indiscrete we know that
 $Bu'$ is a covering by the first paragraph.  Thus, the pullback of
 $Bu$ along the surjective covering map $Br$ is a covering, and it
 follows easily that $Bu$ is a covering.

 Finally, if $H$ is disconnected we just look at one component at a
 time.
\end{proof}

\begin{proposition}\label{prop-coverings-functors}
 Fix a groupoid $H$.  Then the category of coverings $q\:G\xra{}H$ is
 equivalent to the category of functors $X\:H\xra{}\Sets$, and thus
 (by~\cite[Section 1]{qu:haki}) to the category of covering spaces of
 $BH$. 
\end{proposition}
\begin{proof}
 This is a simple translation of Quillen's analysis of coverings of
 $BG$. 

 Suppose we start with a functor $X\:H\xra{}\Sets$.  We then define a
 category $G$ whose objects are pairs $(b,x)$ with $b\in H$ and 
 $x\in X_b$; the morphisms $(b,x)\xra{}(b',x')$ are the maps
 $h\:b\xra{}b'$ in $H$ such that $X_h\:X_b\xra{}X_{b'}$ sends $x$ to
 $x'$.  There is an evident forgetful functor $q\:G\xra{}H$ sending
 $(b,x)$ to $b$; one checks that this is a covering.

 Conversely, suppose we start with a covering $q\:G\xra{}H$.  For each
 $b\in H$, we define $X_b=q^{-1}\{b\}\sse\obj(G)$.  Given a morphism
 $h\:b\xra{}b'$ in $H$ and an element $a\in X_b$, the definition of a
 covering gives a unique morphism $g\:a\xra{}a'$ in $G$ with $q(g)=h$;
 we define a map $X_h\:X_b\xra{}X_{b'}$ by $X_b(a)=a'$.  

 We leave it to the reader to check that these constructions give the
 claimed equivalence.
\end{proof}

We next let $C$ be the class of all coverings, and let $E$ be the
class of functors that are full and essentially surjective.
\begin{proposition}\label{prop-factorisation}
 The pair $(C,E)$ is a factorisation system in $\CG$; in other words
 \begin{itemize}
  \item[(a)] Both $C$ and $E$ contain all identity functors and are 
   closed under composition by isomorphisms on either side.
  \item[(b)] Every functor $u\:G\xra{}H$ can be factored as $u=pr$
   with $p\in C$ and $r\in E$.
  \item[(c)] Every functor in $E$ has the unique left lifting property
   relative to every functor in $C$.  In other words, given functors
   $u$, $w$, $r\in E$ and $p\in C$ making the diagram below commute,
   there is a unique functor $v$ such that $pv=w$ and $vr=u$.
   \begin{diag}
    \node{L} \arrow{s,l}{r} \arrow{e,t}{u} \node{G} \arrow{s,r}{p}\\
    \node{K} \arrow{ne,t,..}{v} \arrow{e,b}{w} \node{H.}
   \end{diag}
 \end{itemize}
\end{proposition}
See~\cite[Exercises 5.5]{bawe:ttt} (for example) for generalities
about factorisation systems.
\begin{proof}
 (a): This is clear.

 (b): Let $u\:G\xra{}H$ be a functor.  We define a new groupoid $K$ as
  follows.  The objects are equivalence classes of triples $(a,b,h)$,
  where $a\in G$ and $b\in H$ and $h\:u(a)\xra{}b$; the equivalence
  relation identifies $(a,b,h)$ with $(a',b',h')$ if and only if
  $b=b'$ and there is a map $g\:a\xra{}a'$ such that
  $h=h'\circ u(g)$.  The maps from $[a,b,h]$ to $[a',b',h']$ are the
  maps $k\:b\xra{}b'$ in $H$ such that there exists a map
  $j\:a\xra{}a'$ in $G$ with $k\circ h=h'\circ u(j)$.  Equivalently,
  $k$ gives a map $[a,b,h]\xra{}[a',b',h']$ if and only if
  $[a',b',h']=[a,b',kh]$.

  There is an evident functor $r\:G\xra{}K$ defined by
  $r(a)=[a,u(a),1_{u(a)}]$.  Given $c=[a,b,h]\in K$ we find that $h$
  can be thought of as a map $r(a)\xra{}b$ in $K$, so $r$ is
  essentially surjective.  Moreover, we find that $K(r(a),r(a'))$ is
  just the image of $G(a,a')$ in $H(u(a),u(a'))$, and thus that $r$ is
  full.  Thus we have $r\in E$.  

  There is also an evident functor $p\:K\xra{}H$ defined by
  $p[a,b,h]=b$.  It is easy to check that $p$ is a covering and $u=pr$
  as required.  In terms of Proposition~\ref{prop-coverings-functors},
  the covering $p$ corresponds to the functor $X\:H\xra{}\Sets$
  defined by $X_b=\pi_0(u\downarrow b)$.

 (c): Suppose we have a square as in the statement of the
  proposition.  We first define a map $v\:\obj(K)\xra{}\obj(G)$ as
  follows.  Suppose that $c\in\obj(K)$.  As $r$ is essentially
  surjective, we can choose $d\in\obj(L)$ and $k\:r(d)\xra{}c$ in $K$.
  We apply $w$ to get $w(k)\:pu(d)=wr(d)\xra{}w(c)$.  As $p$ is a
  covering, there is a unique pair $(a,g)$ with $a\in\obj(G)$ and
  $g\:u(c)\xra{}a$ such that $p(a)=w(c)$ and $p(g)=w(k)$.  We would
  like to define $v(c)=a$.  To check that this is well-defined,
  consider another $d'\in\obj(L)$ and another $k'\:r(d')\xra{}c$,
  giving rise to a unique pair $(a',g')$.  As $r$ is full there exists
  $l\:d'\xra{}d$ such that $k^{-1}k'=r(l)$ and one checks that
  $(a,g\circ u(l))$ has the defining property of $(a',g')$.  Thus
  $a=a'$ as required.  This means that we have a well-defined map
  $v\:\obj(K)\xra{}\obj(G)$ with $pv=w$ on objects.  It is easy to
  check that $vr=u$ on objects as well.

  Now suppose we have a map $m\:c\xra{}c'$ in $K$.  We can choose maps
  $k\:r(d)\xra{}c$ and $k'\:r(d')\xra{}c'$ with $d,d'\in L$.  By the
  definition of $v$ on objects we have maps $g\:u(d)\xra{}v(c)$ and
  $g'\:u(d')\xra{}v(c')$ such that $p(g)=w(k)$ and $p(g')=w(k')$.  As
  $r$ is full we can choose $n\:d\xra{}d'$ such that
  $r(n)=(k')^{-1}mk$.  One then checks that the map
  $g''=g'\circ u(n)\circ g^{-1}\:v(c)\xra{}v(c')$ has $p(g'')=w(m)$.
  As $p$ is faithful, there is at most one map $v(c)\xra{}v(c')$ with
  this property, so $g''$ is independent of the choices made.  We can
  thus define $v$ on morphisms by $v(m)=g''$, so that $pv=w$.  Using
  the faithfulness of $p$, we check easily that $v$ is a functor and
  that $vr=u$.  Thus $v$ fills in the diagram as required.

  Finally suppose that $v'\:K\xra{}G$ is another functor making the
  diagram commute.  We must check that $v'=v$.  As $p$ is faithful it
  is enough to check this on objects.  Given $c\in\obj(K)$ we choose
  $k\:r(d)\xra{}c$ as before and write $a=v'(c)$ and
  $g=v'(k)\:u(d)\xra{}a$.  We then have $p(g)=pv'(k)=w(k)$, so the
  definition of $v$ gives $v(c)=a=v'(c)$ as required.
\end{proof}
\begin{corollary}
 \begin{itemize}
  \item[(i)]   The factorisation in~(b) is unique up to isomorphism.
  \item[(ii)]  $C\cap E$ is precisely the class of isomorphisms in
   $\CG$. 
  \item[(iii)] $C$ and $E$ are closed under compositions and
   retracts.
  \item[(iv)]  $C$ is closed under pullbacks, and $E$ is closed under
   pushouts.
 \end{itemize}
\end{corollary}
\begin{proof}
 See~\cite[Exercises 5.5]{bawe:ttt}.  Of course, in our case, many of
 these things are immediate from the definitions.
\end{proof}

\begin{proposition}\label{prop-quasi-covering}
 A functor $u\:G\xra{}H$ is a quasi-covering if and only if it is
 faithful.  
\end{proposition}
\begin{proof}
 As equivalences and coverings are faithful, we see that
 quasi-coverings are faithful.

 For the converse, let $u\:G\xra{}H$ be faithful.  We can factor $u$
 as $pr$ where $p$ is a covering and $r$ is full and essentially
 surjective, as in Proposition~\ref{prop-factorisation}.  As $u=pr$ is
 faithful we see that $r$ is faithful and thus an equivalence, as
 required. 
\end{proof}

\begin{lemma}\label{lem-factor-backwards}
 Suppose we have functors $L\xra{p}K\xra{v}H$ such that $p$ is a
 covering and $v$ is an equivalence.  Then there is a pullback square
 as follows, in which $q$ is a covering and $u$ is an equivalence.
 \begin{diag}
  \node{L} \arrow{s,l}{p} \arrow{e,t}{u} \node{G} \arrow{s,r}{q}\\
  \node{K} \arrow{e,b}{v} \node{H.}
 \end{diag}
\end{lemma}
\begin{proof}
 We can factor $vp$ as $qu$ with $q$ a covering and $u$ full and
 essentially surjective.  Now consider the following diagram:
 \begin{diag}
  \node{L}
  \arrow{e,t}{w} 
  \arrow{se,b}{p}
  \node{M}
  \arrow{s,l}{\tq}
  \arrow{e,t}{\tv}
  \node{G}
  \arrow{s,r}{q} \\
  \node[2]{K}
  \arrow{e,b}{v}
  \node{H.}
 \end{diag}
 The square is defined to be the pullback of $v$ and $q$, and $w$ is
 the unique functor such that $\tq w=p$ and $\tv w=u$.  By
 Proposition~\ref{prop-coproper} we know that $\tv$ is an equivalence.
 It will thus be enough to show that $w$ is an isomorphism in $\CG$.

 As $u=\tv w$ is full and essentially surjective, and $\tv$ is an
 equivalence, we see that $w$ is full and essentially surjective.
 
 We next show that $w$ is surjective on objects.  Suppose $e\in M$,
 and put $a=\tv(e)\in G$ and $c=\tq(e)\in K$ so that $q(a)=v(c)=b$
 say.  As $u$ is essentially surjective, we can choose $d\in L$ and
 $g\:u(d)\xra{}a$ in $G$.  Thus $q(g)\:vp(d)=qu(d)\xra{}q(a)=v(c)$ in
 $H$.  As $v$ is an equivalence, there is a unique $k\:p(d)\xra{}c$
 such that $v(k)=q(g)$.  As $p$ is a covering, there is a unique
 $d'\in L$ and $l\:d\xra{}d'$ such that $p(d')=c$ and $p(l)=k$.  Thus
 $u(l)g^{-1}\:a\xra{}u(d')$ satisfies
 $q(u(l)g^{-1})=vp(l)v(k)^{-1}=1$.  As $q$ is a covering, it reflects
 identity maps, so $u(l)g^{-1}=1_a$ and $a=u(d')$.  Thus $u$ is
 surjective on objects, as claimed.

 Now consider the following diagram:
 \begin{diag}
  \node{L} \arrow{e,t}{1} \arrow{s,l}{w} \node{L} \arrow{s,r}{p} \\
  \node{M} \arrow{e,b}{\tq} \arrow{ne,b,..}{z} \node{K.}
 \end{diag}
 It follows from Proposition~\ref{prop-factorisation} that there
 is a unique map $z$ making everything commute.  In particular, we
 have $zw=1$.  It follows that $(wz)w=w$ and $w$ is full and
 surjective on objects so $wz=1$.  Thus $w$ is an isomorphism, as
 required. 
\end{proof}

\begin{lemma}\label{lem-cart-quasi}
 Suppose that we have a homotopy cartesian square as follows, in which
 $p$ is a quasicovering.  
 \begin{diag}
  \node{L} \arrow{e,t}{u} \arrow{s,l}{q} \node{G} \arrow{s,r}{p} \\
  \node{K} \arrow{e,b}{v}                \node{H}
 \end{diag}
 Then there is a diagram as follows, in which $p'$ and $q'$ are
 coverings, $r$ and $s$ are equivalences, the bottom square is
 cartesian, the top square commutes up to homotopy, and $p=p'r$,
 $q=q's$. 
 \begin{diag}
  \node{L}  \arrow{e,t}{u} \arrow{s,l}{s}  \node{G}  \arrow{s,r}{r}\\
  \node{L'} \arrow{e,t}{u'}\arrow{s,l}{q'} \node{G'} \arrow{s,r}{p'}\\
  \node{K}  \arrow{e,b}{v}                 \node{H}
 \end{diag}
\end{lemma}
\begin{proof}
 Using Lemma~\ref{lem-homotopy-cartesian}, it is not hard to reduce to
 the case in which $L$ is the standard homotopy pullback of $p$ and
 $v$.  As $p$ is a quasicovering we can factor it as $p=p'r$ where
 $p'$ is a covering and $r$ is an equivalence.  We can then define
 $L'$, $u'$ and $q'$ so that the bottom square is cartesian, which
 implies that $q'$ is a covering.

 Our next task is to define the functor $s$.  An object $d\in L$ is a
 triple $(a,c,h\:p(a)\xra{}v(c))$.  As $p'$ is a covering and
 $h\:p'r(a)\xra{}v(c)$, we see that there is a unique morphism
 $g'\:r(a)\xra{}a'$ in $G'$ such that $p'(a')=v(c)$ and $p'(g')=h$.
 Thus $(a',c)\in L'$ and we can define $s$ on objects by
 $s(d)=(a',c)$.  Note that $u(d)=a$ and $u's(d)=u'(a',c)=a'$ so we can
 define
 \[ \al_d\eqdef g'\: ru(d) \xra{} u's(d). \]

 Next, consider a morphism $(g,k)\:d_0\xra{}d_1$ in $L$, where
 $d_i=(a_i,c_i,h_i\:p(a_i)\xra{}v(c_i))$ for $i=0,1$.  We define
 $a'_i$ and $g'_i$ as above, and define 
 \[ \gb=g'_1\circ r(g)\circ(g'_0)^{-1} \: a'_0 \xra{}a'_1, \]
 so that the following diagram commutes.
 \begin{diag}
  \node{p'r(a_0)} 
  \arrow[3]{e,t}{p'r(g)} 
  \arrow{se,t}{h_0}
  \arrow[2]{s,l}{p'(g'_0)}
  \node[3]{p'r(a_1)}
  \arrow{sw,t}{h_1}
  \arrow[2]{s,r}{p'(g'_1)} \\
  \node[2]{v(c_0)}
  \arrow{e,b}{v(k)}
  \arrow{sw,=}
  \node{v(c_1)}
  \arrow{se,=} \\
  \node{p'(a'_0)}
  \arrow[3]{e,b}{p'(\gb)}
  \node[3]{p'(a'_1).}
 \end{diag}
 We now define $s$ on morphisms by putting $s(g,k)=(\gb,k)$.  It is
 easy to check that this makes $s$ into a functor, and that
 $\al\:ru\xra{}u's$ is a natural map.  Thus, the top square in our
 diagram commutes up to homotopy.  It is also clear that $q's=q$.

 Thus, all that is left is to check that $s$ is an equivalence.  Let
 $d_0$ and $d_1$ be as above, and suppose given $k\:c_0\xra{}c_1$.  As
 $p'$ is faithful and $r$ is an equivalence, we see that there is at
 most one map $g$ making the upper trapezium of the above diagram
 commute, and at most one map $\gb$ making the lower trapezium
 commute.  Moreover, $g$ exists if and only if $\gb$ does, and they
 determine each other by $\gb=g'_1\circ r(g)\circ(g'_0)^{-1}$ and
 $g=r^{-1}((g'_1)^{-1}\circ g\circ g'_0)$.  Note also that
 $L(d_0,d_1)$ is the set of pairs $(g,k)$ such that the top trapezium
 commutes, and $L'(s(d_0),s(d_1))$ is the set of pairs $(\gb,k)$
 making the bottom trapezium commute.  It follows easily that $s$ is
 full and faithful.  

 Now consider an object $d'=(a',c)\in L'$, so $v(c)=p'(a')$.  As $r$
 is essentially surjective we can choose an object $a\in G$ and a map
 $g'\:r(a)\xra{}a'$ in $G'$.  We thus have an object
 $d=(a,c,p'(g')\:p(a)=p'r(a)\xra{}v(c))$ of $L$.  Clearly $s(d)=d'$ so
 $s$ is surjective on objects, and thus an equivalence as claimed.
\end{proof}

\subsection{Cartesian closure}

Let $G$ and $H$ be groupoids, and let $[G,H]$ denote the category of
functors from $G$ to $H$.  It is easy to see that this is a groupoid
and that functors $K\xra{}[G,H]$ biject naturally with functors 
$K\tm G\xra{}H$.  It follows that $\CG$ is cartesian-closed.  One can
also check that this descends to $\CGb$ in the obvious way.  

We next want to check how this works out in the equivalent category
$\CBb$. 

\begin{lemma}
 Suppose that $X$ and $Y$ are objects of $\CB$ and they have the
 homotopy type of CW complexes.  Then the space $C(X,Y)$ of maps from
 $X$ to $Y$ also lies in $\CB$.
\end{lemma}
\begin{proof}
 By well-known results of Milnor, the space $C(S^k,Y)$ also has the
 homotopy type of a CW complex.  Evaluation at the basepoint of $S^2$
 gives a surjective Hurewicz fibration $\ep_Y\:C(S^2,Y)\xra{}Y$ whose
 fibres have the form $\ep_Y^{-1}\{y\}=\Om^2(Y,y)$.  As $Y\in\CB$ we
 know that these fibres are acyclic and so $\ep_Y$ is a weak
 equivalence, and thus a homotopy equivalence.  By a standard result
 (the dual of~\cite[Corollary II.1.12]{ba:ah}, for example) we deduce
 that $\ep_Y$ is fibre-homotopy equivalent to $1_Y$.  One can also see
 that for any $f\:X\xra{}Y$ we have
 \[ \Om^2(C(X,Y),f)\simeq\{g\:X\xra{}C(S^2,Y)\st\ep_Y\circ g=f\} \]
 and our fibre-homotopy equivalence shows that this is contractible.
 The result follows.
\end{proof}

\begin{proposition}\label{prop-B-cartesian}
 If $G,H\in\CG$ then $B[G,H]\simeq C(BG,BH)$ in $\CBb$.
\end{proposition}
\begin{proof}
 It follows from the lemma that $C(BG,BH)\in\CB$.  Recall that
 $B\:\CGb\xra{}\CBb$ is an equivalence.  Thus, for any $K$ we have
 \begin{align*}
  \CBb(BK,C(BG,BH)) &= \CBb(BK\tm BG,BH)                \\
                    &= \CBb(B(K\tm G),BH)               \\
                    &= \CGb(K\tm G,H)                   \\
                    &= \CGb(K,[G,H])                    \\
                    &= \CBb(BK,B[G,H]).
 \end{align*}
 As $B$ is an equivalence we conclude that
 $\CBb(X,C(BG,BH))=\CBb(X,B[G,H])$ for all $X\in\CBb$, and it follows
 by Yoneda's lemma that $C(BG,BH)\simeq B[G,H]$ as claimed.
\end{proof}

\section{Transfers}
\label{sec-transfers}

Let $u\:G\xra{}H$ be a covering with finite fibres.  Then
$Bu\:BG\xra{}BH$ is a finite covering map of spaces, so it is
well-known how to define an associated transfer map 
$Tu\:\Sgi BH_+\xra{}\Sgi BG_+$ of spectra.    
This construction is contravariantly functorial and it converts
disjoint unions to wedges and cartesian products to smash products.
If $p$ is a homeomorphism then $Tp=\Sgi p^{-1}_+$.  The well-known
Mackey property of transfers says that if we have a pullback square as
shown on the left, in which $p$ is a covering, then $q$ is also a
covering and the square on the right commutes up to homotopy.
\begin{diag}
 \node{L}
 \arrow{e,t}{u}
 \arrow{s,l}{q}
 \node{G}
 \arrow{s,r}{p}
 \node[2]{\Sgi BL_+}
 \arrow{e,t}{\Sgi Bu_+}
 \node{\Sgi BG_+}        \\
 \node{K}
 \arrow{e,b}{v}
 \node{H}
 \node[2]{\Sgi BK_+}
 \arrow{n,l}{Tq}
 \arrow{e,b}{\Sgi Bv_+}
 \node{\Sgi BH_+}
 \arrow{n,r}{Tp}
\end{diag}

It will be convenient to extend this to quasicoverings rather than
just coverings.  If $u\:G\xra{}H$ is a quasicovering then we can
factor $u$ as $G\xra{v}K\xra{p}H$ where $v$ is an equivalence and $p$
is a covering.  We then define
\[ Tu=(\Sgi Bv^{-1}_+)\circ Tp. \]
To see that this is well-defined, note (using
Proposition~\ref{prop-factorisation}) that any other such
factorisation has the form $G\xra{wv}L\xra{pw^{-1}}H$ for some
isomorphism $w\:K\xra{}L$.  Using this and the equation 
$Tw=\Sgi Bw_+^{-1}$ we see that 
$(\Sgi B(wv)^{-1}_+)\circ T(pw^{-1})=(\Sgi Bv_+^{-1})\circ Tp$
as required.

Now suppose we have quasicoverings $G\xra{u}H\xra{v}K$; we want to
check that $T(vu)=T(u)T(v)$.  It is easy to reduce to the case where
we have functors $L\xra{p}K\xra{v}H$ such that $p$ is a covering and
$v$ is an equivalence; we need to check that
$T(p)T(v)=T(p)v^{-1}=T(vp)$, where we allow ourselves to write $v$
instead of $\Sgi Bv_+$.  Lemma~\ref{lem-factor-backwards} gives
us a pullback diagram as follows, in which $q$ is a covering and $u$
is an equivalence.
\begin{diag}
 \node{L} \arrow{s,l}{p} \arrow{e,t}{u} \node{G} \arrow{s,r}{q}\\
 \node{K} \arrow{e,b}{v} \node{H.}
\end{diag}
By definition we have $T(vp)=u^{-1}T(q)$.  The Mackey property gives
$u T(p)=T(q)v$ so 
\[ T(p)T(v) = T(p)v^{-1} = u^{-1}T(q) = T(vp) \]
as required.

It is easy to check that in this greater generality we still have
$T(p\amalg q)=T(p)\Wedge T(q)$ and $T(p\tm q)=T(p)\Smash T(q)$.  We
also have an extended Mackey property: if the square on the left is
homotopy-cartesian and $p$ is a quasicovering then $q$ is also a
quasicovering and the right hand square commutes up to homotopy (this
follows easily from Lemma~\ref{lem-cart-quasi}).
\begin{diag}
 \node{L}
 \arrow{e,t}{u}
 \arrow{s,l}{q}
 \node{G}
 \arrow{s,r}{p}
 \node[2]{\Sgi BL_+}
 \arrow{e,t}{\Sgi Bu_+}
 \node{\Sgi BG_+}        \\
 \node{K}
 \arrow{e,b}{v}
 \node{H}
 \node[2]{\Sgi BK_+}
 \arrow{n,l}{Tq}
 \arrow{e,b}{\Sgi Bv_+}
 \node{\Sgi BH_+}
 \arrow{n,r}{Tp}
\end{diag}

\section{The $K(n)$-local category}
\label{sec-Kn-local}

Fix a prime $p$ and an integer $n>0$, and let $K=K(n)$ denote the
associated Morava $K$-theory spectrum.  Let $\CK$ denote the category
of spectra that are local with respect to $K(n)$ in the sense of
Bousfield~\cite{bo:lspe,ra:lrc}, and let $\DD$ be the full subcategory
of strongly dualisable objects in $\CK$.  These categories are studied
in detail in~\cite{host:mkl}.   We write $X\Smash Y$ for the
$K(n)$-localised smash product, which makes $\CK$ into a symmetric
monoidal category.  The unit object is $S:=L_K\Sgi S^0$.  

\begin{definition}
 We define a functor $L\:\CGb\xra{}\DD$ by $LG:=L_K\Sgi BG_+$.  (We
 know from~\cite[Corollary 8.7]{host:mkl} that $LG$ is always
 dualisable, so this lands in $\DD$ as indicated.)  It is clear that
 $L(G\tm H)=LG\Smash LH$ and $L(G\amalg H)=LG\Wedge LH$.
\end{definition}

\begin{definition}
 Let $1$ denote the terminal groupoid (with one object and one
 morphism), and write $\ep$ for the unique functor $G\xra{}1$.  Let
 $\dl\:G\xra{}G\tm G$ be the diagonal functor.  Define
 \[ b_G = (LG\Smash LG \xra{L_KT\dl} LG \xra{L\ep} S). \]
 It is not hard to see that $b_{G\tm H}=b_G\Smash b_H$ and
 $b_{G\amalg H}=b_G\oplus b_H$.
\end{definition}

The following result is the key to the whole paper.
\begin{proposition}
 For any finite groupoid $G$, the map $b_G$ is an inner product on
 $LG$.
\end{proposition}
\begin{proof}
 We can easily reduce to the case where $G$ is a group rather than a
 groupoid.  It was observed in the proof of of~\cite[Corollary
 8.7]{host:mkl} that a certain map $c_G\:LG\xra{}DLG$ (arising from the
 Greenlees-May theory of generalised Tate spectra) is an isomorphism.
 It is thus enough to show that $c_G=b_G^\#$.  

 We will need some notation.  Firstly, we will need to consider
 various unlocalised spectra, so in this proof only we write $S$ for
 the ordinary, unlocalised sphere spectrum, and $\hS$ for $L_{K(n)}S$.
 Similarly, we write $X\Smash Y$ for the unlocalised smash product and
 $X\hSmash Y=L_{K(n)}(X\Smash Y)$.  Next, we will work partially in
 the equivariant categories of $G$-spectra and $G^2$-spectra, indexed
 over complete universes~\cite{lemast:esh}.  We write $S'$ and $S''$
 for the corresponding $0$-sphere objects.  Also, we can regard $\hS$
 as a naive $G$-spectrum with trivial action and then extend the
 universe to obtain a genuine $G$-spectrum, which we denote by
 $\hS'$.  We define a genuine $G^2$-spectrum $\hS''$ in the 
 analogous way.

 We next recall the definition of $c_G$.  It is obtained from a
 certain map $d_G\:\hS\Smash BG_+\xra{}F(BG_+,\hS)\simeq DLG$ by
 observing that $DLG$ is $K(n)$-local and that any map from
 $\hS\Smash BG_+$ to a $K(n)$-local spectrum factors uniquely through
 $L_{K(n)}(\hS\Smash BG_+)\simeq LG$.  It will be enough to check that
 $d_G$ is adjoint to $1\Smash b'_G\:\hS\Smash BG^2_+\xra{}\hS$, where
 $b'_G$ is the composite
 \[ \Sgi BG^2_+ \xra{\psi^!} \Sgi BG_+ \xra{B\ep} S. \]

 We thus need to show that two elements of the group
 $[\hS\Smash BG^2_+,\hS]$ are equal.  Theorem~II.4.5
 of~\cite{lemast:esh} (applied to $G^2/G^2\simeq 1$) gives a natural 
 isomorphism
 \[ [\hS\Smash BG^2_+,\hS] \simeq [\hS''\Smash EG^2_+,\hS'']_{G^2}.
 \]
 Let $\zt\:EG_+\xra{}S^0$ and $\xi\:G^2/\Dl_+\xra{}S^0$ be the
 collapse maps.  Desuspending Construction~II.5.1 of~\cite{lemast:esh}
 gives a pretransfer map $t\:S''\xra{}\Sgi G^2/\Dl_+$ of genuine
 $G^2$-spectra.  By smashing this with $EG^2_+$ and passing to orbits
 we get the transfer map $\psi^!\:\Sgi BG^2_+\xra{}\Sgi B\Dl_+$.  Using
 this and the proof of~\cite[Theorem II.4.5]{lemast:esh} we find that 
 $1\Smash b'_G$ corresponds to the composite
 \[ \hS''\Smash EG^2_+ \xra{1\Smash\zt\Smash\zt\Smash t}
    \hS''\Smash G^2/\Dl_+ \xra{1\Smash\xi} \hS''
 \]
 in $[\hS''\Smash EG^2_+,\hS'']_{G^2}$.  

 We now return to the definition of $d_G$.   We have a map
 \[ \hS'\Smash EG_+\xra{1\Smash \zt}\hS'\xra{\zt^*}F(EG_+,\hS') \]
 of $G$-spectra.  We next apply the fixed point functor, noting that
 $F(EG_+,\hS')^G=F(BG_+,\hS)$ and that~\cite[Theorem
 II.7.1]{lemast:esh} gives an equivalence $\ttt\:\hS\Smash
 BG_+\xra{}(\hS'\Smash EG_+)^G$.  The resulting map
 $\hS\Smash BG_+\xra{}F(BG_+,\hS)$ is $d_G$ (see~\cite[Section
 5]{grma:gtc}).  To understand this better, we need to follow through
 the construction of $\ttt$.  We use the notation of~\cite[Section
 II.7]{lemast:esh}, noting that in our case we have $N=G$.  The
 construction uses the group $\Gm=G\tm_cN$, the semidirect product of
 $G$ with $N$ using the action by conjugation.  There are two natural
 maps $\ep,\tht\:\Gm\xra{}G$ given by $\ep(g,n)=g$ and $\tht(g,n)=gn$.
 In our case we find that the resulting map $\Gm\xra{}G^2$ is an
 isomorphism, so we can replace $\Gm$ by $G^2$ everywhere.  The
 subgroup $\Pi$ becomes $1\tm G$, the standard embedded copy $G\tm_c1$
 of $G$ becomes $\Dl$, and the maps $\ep$ and $\tht$ become the
 projections $\pi_0,\pi_1\:G^2\xra{}G$.  The relevant spectrum $D$ is
 $\hS\Smash EG_+$, so $i_*\tht^*D=\hS'\Smash\pi_1^*EG_+$ and
 $j_*i_*\tht^*D=\hS''\Smash\pi_1^*EG_+$.  The map $\ttt$ is obtained
 from 
 \[ 1\Smash t\:\hS''\Smash\pi_1^*EG_+\xra{}
               \hS''\Smash\pi_1^*EG_+\Smash G^2/\Dl_+
 \]
 by shrinking the universe, passing to orbits and adjointing as
 described in~\cite[Construction II.7.5]{lemast:esh}.  It follows that
 $d_G$ is obtained from the composite
 \[ \hS''\Smash\pi_1^*EG_+\xra{1\Smash t}
    \hS''\Smash\pi_1^*EG_+\Smash G^2/\Dl_+
    \xra{1\Smash\zt\Smash 1\Smash\zt^*}
    \hS''\Smash G^2/\Dl_+\Smash F(\pi_0^*EG_+,S)
 \]
 by a similar procedure.  We can identify $EG^2_+$ with
 $\pi_0^*EG_+\Smash\pi_1^*EG_+$, and we find that the adjoint of
 $d_G$ is obtained by applying another similar procedure to the
 map
 \[ 
  \hS''\Smash EG^2_+\xra{1\Smash\zt\Smash\zt\Smash t}
  \hS''\Smash G^2/\Dl_+.
 \]
 This procedure amounts to just composing with
 $\xi\:G^2/\Dl_+\xra{}S^0$ and using our isomorphism
 $[\hS\Smash BG^2_+,\hS]\simeq[\hS''\Smash EG^2_+,\hS'']_{G^2}$.  It
 follows that the adjoint of $d_G$ is $b'_G$, as required.
\end{proof}

\begin{definition}
 For any functor $u\:G\xra{}H$ we put $Ru=(Lu)^t\:LH\xra{}LG$.  
\end{definition}

\begin{proposition}\label{prop-transfer-adjoint}
 If $u\:G\xra{}H$ is faithful then $Ru=L_KTu$.
\end{proposition}
\begin{proof}
 Let $\psi_H\:H\xra{}H\tm H$ be the diagonal map.  We first claim that
 the following square is homotopy-cartesian:
 \begin{diag}
  \node{G} \arrow{s,l}{(1,u)} \arrow{e,t}{u}
  \node{H} \arrow{s,r}{\psi_H}                   \\
  \node{G\tm H} \arrow{e,b}{u\tm 1} \node{H\tm H.}
 \end{diag}
 To see this, let $K$ be the homotopy pullback of the functors $\psi_H$
 and $u\tm 1$.  The square is clearly cartesian, which means that $G$
 embeds as a full subcategory of $K$; we need only check that the
 inclusion is essentially surjective.  The objects of $K$ are
 $5$-tuples $(a,b,c,h,k)$ where $a\in G$ and $b,c\in H$ and
 $h\:u(a)\xra{}c$ and $k\:b\xra{}c$.  The maps from $(a,b,c,h,k)$ to
 $(a',b',c',h',k')$ are triples $(r,s,t)$ making the following diagram
 commute: 
 \begin{diag}
  \node{u(a)}  \arrow{e,t}{h}  \arrow{s,l}{u(r)} 
  \node{c}                     \arrow{s,l}{s}
  \node{b}     \arrow{w,t}{k}  \arrow{s,r}{t}   \\
  \node{u(a')} \arrow{e,b}{h'}
  \node{c'}
  \node{b'.}   \arrow{w,b}{k'}
 \end{diag}
 The canonical functor $v\:G\xra{}K$ is given by
 $v(a)=(a,u(a),u(a),1,1)$.  We define $w\:K\xra{}G$ by
 $w(a,b,c,h,k)=a$.  Then $wv=1$, and we have a natural map
 $vw(a,b,c,h,k)\xra{}(a,b,c,h,k)$ given by $(1,k,k^{-1}h)$.  This
 proves that $v$ is an equivalence, and if we compose it with the
 projections $K\xra{}G\tm H$ and $K\xra{}H$ we get the functors
 $(1,u)$ and $u$.  This proves that our original square is
 homotopy-cartesian, so the Mackey property tells us that
 \[ T(\psi_H)\circ(Bu\tm 1)=Bu\circ T(1,u) \:
     \Sgi B(G\tm H)_+\xra{}\Sgi BH_+.
 \]
 We now use the fact that $(1,u)=(1\tm u)\circ\psi_G$ and compose with
 the projection $\Sgi BH_+\xra{}S^0$ to get
 \[ b_H\circ(Bu\Smash 1)=\ep_H\circ Bu\circ(T\psi_G)\circ(1\Smash Tu).
 \]
 We next note that $\ep_H\circ Bu=\ep_G$ and $K(n)$-localise to
 conclude that $b_H\circ(Lu\Smash 1)=b_G\circ(1\Smash L_KTu)$, as
 required. 
\end{proof}

We can thus think of the maps $Ru$ as generalised transfers.  It turns
out that we also have a generalised Mackey property.
\begin{proposition}\label{prop-mackey-general}
 If we have a homotopy-cartesian square as shown on the left, then the
 diagram on the right commutes.
 \begin{diag}
  \node{M}
  \arrow{e,t}{u}
  \arrow{s,l}{t}
  \node{G}
  \arrow{s,r}{s}
  \node[2]{LM}
  \arrow{e,t}{Lu}
  \node{LG}        \\
  \node{K}
  \arrow{e,b}{v}
  \node{H}
  \node[2]{LK}
  \arrow{n,l}{Rt}
  \arrow{e,b}{Lv}
  \node{LH}
  \arrow{n,r}{Rs}
 \end{diag}
\end{proposition}
\begin{proof}
 We may assume that the square is actually a pullback square of
 fibrations (see Remark~\ref{rem-replace-square}), so in particular it
 commutes on the nose.  As $b_H$ is a perfect pairing, it suffices to
 check that $b_G\circ(1\Smash (Rs)(Lv))=b_G\circ(1\Smash (Lu)(Rt))$.
 By transposition, this is equivalent to
 \[ b_H\circ(Ls\Smash Lv)=b_M\circ(Ru\Smash Rt)\:LG\Smash LK\xra{}S.
 \] 
 To verify this, we consider the following diagram:
 \begin{diag}
  \node{M} \arrow{e,t}{vt=su} \arrow{s,l}{(u,t)}
  \node{H} \arrow{s,r}{\psi_H} \\
  \node{G\tm K} \arrow{e,b}{s\tm v}
  \node{H\tm H.}
 \end{diag}
 We claim that this is homotopy-cartesian.  It is clearly cartesian,
 so it suffices (as in the previous proof) to show that the obvious
 functor from $M$ to the homotopy pullback is essentially surjective.
 Suppose we are given an object of the homotopy pullback, in other
 words a $5$-tuple $d=(a,b,c,k,l)$ where $a\in G$, $b\in H$, $c\in K$
 and $s(a)\xra{k}b\xla{l}v(c)$.  As $s$ is a fibration we can choose
 $a'\in G$ and $g\:a\xra{}a'$ such that $s(a')=v(c)$ and
 $s(g)=l^{-1}k$.  Thus $d'\eqdef(a',s(a')=v(c),c,1,1)$ is an object of
 $M$ and the following diagram gives an isomorphism $d\xra{}d'$:
 \begin{diag}
  \node{s(a)} \arrow{e,t}{k} \arrow{s,l}{s(g)} 
  \node{b}                   \arrow{s,l}{l^{-1}}
  \node{v(c)} \arrow{w,t}{l} \arrow{s,r}{v(1)} \\
  \node{s(a')} \arrow{e,=} \node{v(c)} \node{v(c)} \arrow{w,=}
 \end{diag}
 This shows that our square is homotopy-cartesian.  The vertical
 functors are faithful and thus are quasicoverings, so the Mackey
 property tells us that 
 \[ (R\psi_H)\circ(Ls\Smash Lv)=L(su)\circ R(u,t)\:LG\Smash LK\xra{}LH.
 \]
 We next compose with the map $L\ep_H\:LH\xra{}S$, noting that
 $\ep_Hsu=\ep_M\:M\xra{}1$ and that
 $(u,t)=(u\tm t)\psi_M\:M\xra{}G\tm K$.  We conclude that
 $b_H\circ(Ls\Smash Lv)=b_M\circ(Ru\Smash Rt)$, as required.
\end{proof}

\begin{theorem}\label{thm-groupoids-frobenius}
 For any finite groupoid $G$, the maps $(R\ep,R\psi,L\ep,L\psi)$ make
 $LG$ into a Frobenius object.  
\end{theorem}
\begin{proof}
 This is formally identical to the proof of
 Proposition~\ref{prop-manifolds-frobenius}; we need only check that
 the following square is homotopy-cartesian, and that is easy.
 \begin{diag}
  \node{G}           \arrow{e,t}{\psi}      \arrow{s,l}{\psi}
  \node{G\tm G}                            \arrow{s,r}{1\tm\psi}\\
  \node{G\tm G}      \arrow{e,b}{\psi\tm 1} 
  \node{G\tm G\tm G.}
 \end{diag}
 Alternatively, the result can be deduced from the proof of
 Proposition~\ref{prop-psi}.
\end{proof}

\begin{definition}
 Given a finite groupoid $G$, define $\Lm G=[\Z,G]$;
 Proposition~\ref{prop-B-cartesian} tells us that $B\Lm G$ is
 homotopy-equivalent to the free loop space on $BG$.  The objects of
 $\Lm G$ are pairs $(a,u)$ where $u\in G(a,a)$, and the maps from
 $(a,u)$ to $(b,v)$ are maps $g\:a\xra{}b$ such that $v=gug^{-1}$.  It
 is thus easy to see that $\pi(a,u)=a$ gives a functor $\Lm G\xra{}G$,
 and that this is actually a covering.  If $G$ is a group then 
 $\Lm G$ is equivalent the disjoint union of the groups $Z_G(g)$ as
 $g$ runs over the conjugacy classes in $G$, so the free loop space on
 $BG$ is $\coprod BZ_G(g)$; this is actually well-known, and a more
 elementary account appears in~\cite[Section 2.12]{be:rcii}, for
 example. 
\end{definition}
\begin{remark}
 It is important to distinguish between $[\Z,G]$ and $[\Zp,G]$; see
 Section~\ref{sec-warnings} for more discussion of this.
\end{remark}

We can now identify the maps $\tht=\ep\mu\psi\:LG\xra{}S$ and
$\al=\mu\psi\eta\:S\xra{}LG$ discussed in
Proposition~\ref{prop-frobenius-trace}.
\begin{proposition}\label{prop-trace-loops}
 We have $\tht=(L\ep_{\Lm G})(R\pi)$ and $\al=(L\pi)(R\ep_{\Lm G})$.
\end{proposition}
\begin{proof}
 The key point is that the following square is homotopy-cartesian:
 \begin{diag}
  \node{\Lm G} \arrow{e,t}{\pi} \arrow{s,l}{\pi}
  \node{G}                      \arrow{s,r}{\psi} \\
  \node{G}     \arrow{e,b}{\psi} \node{G^2.}
 \end{diag}
 To see this, let $H$ be the homotopy pullback of $\psi$ and $\psi$.
 The objects of $H$ are tuples $(a,b,u,v)$ where $a,b\in G$ and
 $u,v\:a\xra{}b$.  The morphisms from $(a,b,u,v)$ to $(a',b',u',v')$
 are pairs $(g,h)$ where $g\:a\xra{}a'$ and $g\:b\xra{}b'$ and
 $hu=u'g$ and $hv=v'g$.  We can define a functor
 $\phi\:\Lm G\xra{}H$ by $(a,u)\mapsto(a,a,u,1)$ and a functor $\xi$
 in the opposite direction by $(a,b,u,v)\mapsto(a,v^{-1}u)$.  We find
 that these are equivalences and that either projection $H\xra{}G$
 composed with $\phi$ is just $\pi$; it follows that the square is
 homotopy-cartesian, as claimed.  We conclude that
 $\mu\psi=(R\psi)(L\psi)=(L\pi)(R\pi)$.  We also know from
 Proposition~\ref{prop-frobenius-trace} that
 $\tht=\ep\mu\psi=(L\ep)(R\psi)(L\psi)$ and
 $\al=\mu\psi\eta=(R\psi)(L\psi)(R\ep)$.  Everything now follows from the
 evident fact that $\ep_G\pi=\ep_{\Lm G}\:\Lm G\xra{}1$.
\end{proof}

We conclude this section by discussing the case of a finite abelian
group $A$, considered as a groupoid with one object.  There is then a
unique functor $\zt\:1\xra{}A$, and also a division homomorphism
$\nu\:A\tm A\xra{}A$ given by $\nu(g,h)=gh^{-1}$.  
\begin{proposition}\label{prop-bA-abelian}
 We have $b=(R\zt)(L\nu)\:LA\Smash LA\xra{}S$.  We also have
 $\al=|A|\eta$ and $\tht=|A|\ep$.
\end{proposition}
\begin{proof}
 We have a commutative diagram as follows, which is easily seen to be
 both cartesian and homotopy-cartesian:
 \begin{diag}
  \node{A} \arrow{s,l}{\psi} \arrow{e,t}{\ep}\node{1}\arrow{s,r}{\zt}\\
  \node{A\tm A} \arrow{e,b}{\nu} \node{A.}
 \end{diag}
 The vertical functors are faithful and thus are quasicoverings.  The
 Mackey property now tells us that $b=(L\ep)(R\psi)=(R\zt)(L\nu)$ as
 claimed. 

 Next, consider the groupoid $\Lm A=[\Z,A]$.  It is easy to see that
 this is just a disjoint union of $|A|$ copies of $A$, and that the
 functor $\pi\:\Lm A\xra{}A$ just sends each copy isomorphically to
 $A$.  The remaining claims now follow easily from
 Proposition~\ref{prop-trace-loops}. 
\end{proof}

\section{Inner products in cohomology}
\label{sec-inner-cohomology}

We next study $E^*BG$ for suitable cohomology theories $E$.

If $p$ is an odd prime, let $E$ be a $p$-local commutative ring
spectrum such that
\begin{itemize}
 \item[(a)] $E^0$ is a complete local Noetherian ring
 \item[(b)] $E^1=0$
 \item[(c)] $E^2$ contains a unit
 \item[(d)] The associated formal group over $\spec(E^0/\mxi)$ has
  height $n$.
\end{itemize}
Here~(d) makes sense because~(b) and~(c) force the coefficient ring
$E^*$ to be concentrated in even degrees, so the Atiyah-Hirzebruch
spectral sequence for $E^*\cpi$ collapses, so $E$ is automatically
complex-orientable.  In the language of~\cite[Section 2]{grst:vlc},
our assumption is that $E$ is a $K(n)$-local admissible ring spectra.

In the case $p=2$ we would like to allow $E$ to be a two-periodic
version of $K(n)$, but this is not commutative.  We therefore relax
the requirement that $E$ be commutative and assume instead that there
is a derivation $Q\:E\xra{}\Sg E$ and an element $v\in\pi_2E$ such
that $2v=0$ and
\[ ab - ba = vQ(a)Q(b), \]
so that $E$ is quasicommutative in the sense of~\cite[Definition
8.1.1]{st:fsfg}.  This is of course satisfied if $E$ is commutative,
with $Q=0$ and $v=0$.  Other examples, including the two-periodic
version of $K(n)$, can most easily be produced by the methods
of~\cite{st:pmm}, which also contains detailed references to previous
work in this direction.

We consider $LG$ as a Frobenius object just as in the previous
section.  As usual we use the maps
$S\xla{L\ep}LG\xra{L\psi}LG\Smash LG$ to make $E^0LG=E^0BG$ into a ring
and $E^*LG=E^*BG$ into a graded ring.  We also use $(R\ep)\:S\xra{}LG$
to give a map $\ep:=(R\ep)^*\:E^0BG\xra{}E^0$, which in turn gives a
bilinear form $b(x,y)=\ep(xy)$ on $E^0BG$.

\begin{remark}
 If $G$ is a group then the inclusion of the trivial group gives a map
 $\zt\:1\xra{}G$ and thus an augmentation map
 $(L\zt)^*\:E^0BG\xra{}E^0$.  In other contexts this is often denoted
 by $\ep$, but it is \emph{not} the same as the map $\ep$ defined
 above.
\end{remark}

We say that $G$ is \emph{$E$-good} if $E^0LG$ is free of finite rank
over $E^0$ and $E^1LG=0$.  If so then we have a K\"unneth isomorphism
$E^0(LG\Smash LG)=E^0LG\ot_{E^0}E^0LG$.  Using this and
Theorem~\ref{thm-groupoids-frobenius} we find that the above maps make
$E^0LG$ into a Frobenius object in the compact-closed category of
finitely generated free modules over $E^0$.  In particular, we deduce
that our bilinear form is an inner product.  A functor $u\:G\xra{}H$
gives a ring map $u^*\:E^0BH\xra{}E^0BG$ induced by $Lu\:LG\xra{}LH$,
and also a map $u_!:=(Ru)^*\:E^0BG\xra{}E^0BH$ that is adjoint to
$u^*$.  If $u$ is the inclusion of a subgroup in a group then $u_!$ is
the corresponding transfer map (by
Proposition~\ref{prop-transfer-adjoint}).  The adjointness of $u_!$
and $u^*$ is thus a version of Frobenius reciprocity.

As usual we have a trace map $\tht\:E^0BG\xra{}E^0$ (which can be
computed using only the ring structure) and an element
$\al=\mu\psi(1)\in E^0BG$.  Proposition~\ref{prop-trace-loops} tells
us how to compute $\al$ in terms of ordinary transfers, and
Proposition~\ref{prop-frobenius-trace} tells us that
$\ep(x\al)=\tht(x)$.  We will see later that $\al$ becomes invertible
in $\Q\ot E^0BG$, so the previous equation characterises $\ep$ up to
torsion terms.

Now let $A$ be a finite Abelian group.  It is known that such groups
are $E$-good for all $E$; see~\cite[Proposition 2.9]{grst:vlc} for a
proof in the present generality, although the basic idea of the proof
is much older~\cite{la:ccc,hokura:ggc}.  We know from
Proposition~\ref{prop-bA-abelian} that $\al=|A|$ in this context so
that $|A|\ep(x)=\tht(x)$.  We next give another formula for $\ep$ that
is more useful when $p=0$ in $E^0$.  It is easy to see that $\ep_{A\tm
  B}=\ep_A\ot\ep_B$, and if $|A|$ is coprime to $p$ then $E^0BA=E^0$
with $\ep_A=|A|.1\:E^0\xra{}E^0$.  It is thus enough to treat the case
where $A=C_{p^m}$ for some $m>0$.

Choose a complex orientation $x\in\tE^0\cpi$, or equivalently a
coordinate on the associated formal group $\GG$.  This gives a formal
group law $F$ with associated $p^m$-series $[p^m](x)$, and we have
\[ E^0B(C_{p^m})=\O_{\GG(m)}=\fps{E^0}{x}/[p^m](x). \]
There is a unique invariant differential form $\om$ on $G$ that agrees
with $dx$ at zero: if $E^0$ is torsion-free this is most easily
expressed as $\om=d(\log_F(x))=\log'_F(x)dx$.  Given a function
$f\in\O_\GG=\fps{E^0}{x}$ we get a meromorphic differential form
$f\om/[p^m](x)\in\CM\Om^1_{\GG}$, and the residue of this form clearly
only depends on $f$ modulo $[p^m](x)$.  (See~\cite[Sections 5.3 and
5.4]{st:fsfg} for an exposition of meromorphic forms and their
residues.)

\begin{proposition}\label{prop-frobenius-residue}
 The canonical Frobenius form on $E^0BC_{p^m}$ is given by
 $\ep(f)=\res(f\om/[p^m](x))$. 
\end{proposition}
\begin{proof}
 For any $E$-good group $G$, we can define 
 \[ c := \psi\eta(1) = 
    \text{tr}_\Dl^{G^2}(1) \in E^0(BG^2)=E^0BG\ot_{E^0}E^0BG.
 \]
 We see from Scholium~\ref{sch-determine-ep} that
 $\ep\:E^0BG\xra{}E^0$ is the unique map such that
 $(\ep\ot 1)(c)=1\in E^0BG$.

 Now take $G=C_{p^m}$, so $E^0BG=\fps{E^0}{x}/[p^m](x)$ and
 $E^0BG^2=\fps{E^0}{x,y}/([p^m](x),[p^m](y))$.  Write
 $\dps{m}(t)=[p^m](t)/t\in\fps{E^0}{t}$; we know from~\cite[Section
 4]{qu:epc} that $\tr_1^G(1)=\dps{m}(x)$ (a simpler proof appears
 in~\cite{tu:dlo}).  Put $z=x-_Fy$; it follows from
 Proposition~\ref{prop-bA-abelian} that $c=\dps{m}(z)$.

 Now consider the form 
 \[ \gm = c\om/[p^m](x)\in 
     \frac{\fps{E^0}{y}}{[p^m](y)} \ot_{E^0} \CM\Om^1_{\GG},
 \]
 so that $\res(\gm)\in\fps{E^0}{y}/[p^m](y)$.  In view of the above,
 it will suffice to check that $\res(\gm)=1$.

 For this, we note that $[p^m](y)=0$ so $zc=[p^m](z)=[p^m](x)$ so
 $z\gm=\om$ so $\gm=\om/z$.  Now, $\om=g(x)dx$ for some power series
 $g$ with $g(0)=1$ and this differential is invariant under
 translation, which implies that $\om=g(z)dz$ also.  Thus
 $\res(\gm)=\res(g(z)dz/z)=g(0)=1$ as required.
\end{proof}
\begin{corollary}
 Let $E$ be the usual two-periodic version of $K(n)$ (with $n>1$), and
 let $x$ be the usual $p$-typical orientation.  Then the Frobenius
 form on the ring
 \[ E^0BC_{p^m}=\fps{E^0}{x}/x^{p^{nm}}=E^0\{x^k\st 0\leq k<p^{nm}\}
 \] 
 is given by $\ep(x^k)=0$ for $k<p^{nm-1}$ and $\ep(x^{p^{nm}-1})=1$.
 In the case $n=1$ we have $\ep(x^{p^m-p^j})=1$ for $0\leq j\leq m$
 and $\ep(x^k)=0$ for all other $k$.
\end{corollary}
\begin{proof}
 For the integral two-periodic version of $K(n)$ we have
 $\log_F(x)=\sum_{k\geq 0}x^{p^{nk}}/p^k$.  When $n>1$ it follows
 easily that $\om=\log'_F(x)dx=dx\pmod{p}$.  We also have
 $[p^m](x)=x^{p^{mn}}$ so $\ep(x^k)=\res(x^{k-p^{mn}}dx)$ and the
 claim follows easily.  In the case $n=1$ we have
 $\om=\sum_{k\geq 0}x^{p^k-1}dx$ and the stated formula follows in the
 same way.
\end{proof}

\section{Character theory}
\label{sec-characters}

Let $G$ be a finite groupoid.  Write $C(G):=\Q\{\pi_0G\}$ for the
rational vector space freely generated by the set of isomorphism
classes of objects of $G$.  Given $a\in G$ we write $[a]$ for the
corresponding basis element in $C(G)$.  We define a bilinear form on
$C(G)$ by
\[ ([a],[b]) := |G(a,b)|. \]
It is convenient to write $G(a):=G(a,a)$ and to introduce the elements
$[a]':=[a]/|G(a)|$, so that $([a],[a]')=1$.  We also write
$C(G)^*=\Hom_\Q(C(G),\Q)=F(\pi_0G,\Q)$ for the dual of $C(G)$.  Given
a functor $u\:G\xra{}H$ we define $Lu\:C(G)\xra{}C(H)$ by
$(Lu)[a]=[u(a)]$, and we let $Ru\:C(H)\xra{}C(G)$ be the adjoint of
this, so that
\[ (Ru)[b]' = \sum_{[a]\st u(a)\simeq b} [a]'. \]
The sum here is indexed by isomorphism classes of objects $a\in G$
such that $u(a)$ is isomorphic to $b$ in $H$.  

We next show that these constructions have the expected Mackey
property.
\begin{proposition}\label{prop-mackey-characters}
 If we have a homotopy-cartesian square as shown on the left, then the
 diagram on the right commutes.
 \begin{diag}
  \node{M}
  \arrow{e,t}{u}
  \arrow{s,l}{t}
  \node{G}
  \arrow{s,r}{s}
  \node[2]{C(M)}
  \arrow{e,t}{Lu}
  \node{C(G)}        \\
  \node{K}
  \arrow{e,b}{v}
  \node{H}
  \node[2]{C(K)}
  \arrow{n,l}{Rt}
  \arrow{e,b}{Lv}
  \node{C(H)}
  \arrow{n,r}{Rs}
 \end{diag}
\end{proposition}
\begin{proof}
 We may assume that the square is actually a pullback square of
 fibrations (see Remark~\ref{rem-replace-square}), so in particular it
 commutes on the nose.  Fix $c\in K$, so $(Lv)[c]=[vc]=|H(vc)|[vc]'$.
 We need to check that $(Rs)(Lv)[c]=(Lu)(Rt)[c]$.  Because $s$ is a
 fibration, any isomorphism class in $G$ that maps to $[vc]$ in $H$
 has a representative $a\in G$ such that $sa=vc$.  Using this, we find
 that 
 \[ (Rs)(Lv)[c] = \sum_{[a]\st sa=vc}|H(vc)||G(a)|^{-1}[a]. \]
 We also know that $t$ is a fibration, so every isomorphism class in
 $M$ that maps to $c$ contains a representative $d$ with $t(d)=c$, in
 other words $d$ has the form $(a,c)$ for some $a\in G$ with $sa=vc$.
 It follows that
 \[ (Lu)(Rt)[c] = \sum_{[a,c]\st sa=vc}|M((a,c))|^{-1}|K(c)|[a]. \]
 Fix $a\in G$ with $sa=vc$.  The coefficient of $[a]$ in $(Rs)(Lv)[c]$
 is then $|H(vc)||G(a)|^{-1}$.  For $(Lu)(Rt)[c]$ we need to be more
 careful, because there will typically be objects $a'\in G$ with
 $[a']=[a]\in\pi_0G$ but $[a',c]\neq[a,c]\in\pi_0M$.  Put
 \[ X=\{a'\in G\st a'\simeq a\text{ and } sa'=vc\}/\sim, \]
 where $a'\sim a''$ iff there exist $g\:a'\xra{}a''$ and $k\:c\xra{}c$
 such that $sg=vk\:vc\xra{}vc$.  It is easy to see that $a'\sim a''$
 iff $(a',c)\simeq(a'',c)$ in $M$, and it follows that the coefficient
 of $[a]$ in $(Lu)(Rt)[c]$ is
 \[ \lm := \sum_{[a']\in X} |M((a',c))|^{-1}|K(c)|. \]
 To analyse this further, we introduce the set
 \[ Y=\{(a',g')\st a'\in G\;,\;sa'=vc\text{ and }g'\:a'\xra{}a\}/\sim,
 \]
 where $(a',g')\sim(a'',g'')$ iff $sg'=sg''\:vc\xra{}vc$.  Using the
 fact that $s$ is a fibration, one checks that the map
 $[a',g']\mapsto sg'$ gives a bijection $Y\simeq H(vc)$, so that
 $|Y|=|H(vc)|$.  On the other hand, there is an evident projection
 $\pi\:Y\xra{}X$ sending $[a',g']$ to $[a']$.  If $\pi[a'',g'']=[a']$
 then we can choose $f\:a'\xra{}a''$ and $k\:c\xra{}c$ such that
 $sf=vk$, and then observe that there is a unique $g\:a'\xra{}a'$ such
 that $f=(g'')^{-1}g'g$.  One checks that the coset $(g,k).M((a',c))$
 depends only on the equivalence class $[a'',g'']$  and that this
 construction gives a bijection 
 \[ \pi^{-1}\{[a']\}\simeq (G(a')\tm K(c))/M((a',c)). \]
 Note also that $|G(a')|=|G(a)|$ because $a'\simeq a$.  It follows
 that 
 \begin{align*}
  \lm &= \sum_{[a',g']\in Y}
          |G(a')\tm K(c)/M((a',c))|^{-1} |M((a',c))|^{-1} |K(c)|\\
      &= \sum_{[a',g']\in Y} |G(a)|^{-1} \\
      &= |H(vc)||G(a)|^{-1}.
 \end{align*}
 This is the same as the coefficient of $[a]$ in $(Rs)(Lv)[c]$, as
 required. 
\end{proof}

If we let $1\xla{\ep}G\xra{\psi}G^2$ be the obvious functors, then it
follows easily that the maps $(L\ep)^*$, $(L\psi)^*$, $(R\ep)^*$,
$(R\psi)^*$ make $C(G)^*$ into a Frobenius algebra over $\Q$.  The
Frobenius form is just
\[ \ep(f)=(f,1)=\sum_{[a]\in\pi_0G} f(a)/|G(a)|. \]

Next, let $\Lm$ denote the group $(\Qp/\Zp)^n$, whose dual is
$\Lm^*=\Hom(\Lm,\Qp/\Zp)\simeq\Zp^n$.  We regard $\Lm^*$ as a groupoid
with one object in the usual way.  We also consider
$\Lm(m)=\ker(\Lm\xra{p^m}\Lm)$, so $\Lm(m)^*=\Lm^*/p^m$.

Note that if $G$ is a finite group and $\lm\:\Lm^*\xra{}G$ then the
image of $\lm$ is Abelian and $p$-local (because it is a quotient of
$\Lm^*$) and finite (because it is a subgroup of $G$).  It follows
that $\lm(p^m\Lm^*)=1$ for large $m$, so that $\lm$ is automatically
continuous.  It follows that $\Hom(\Lm^*,G)$ bijects with the set of
$n$-tuples of commuting elements of $G$ of $p$-power order.  More
generally, if $G$ is a finite groupoid then a functor
$\lm\:\Lm^*\xra{}G$ consists of an object $a\in G$ together with an
$n$-tuple of commuting $p$-elements of the group $G(a,a)$.

The generalised character theory of Hopkins, Kuhn and
Ravenel~\cite{hokura:ggc,hokura:mkc} can be repackaged slightly to
relate $\Q\ot E^0BG$ to $C([\Lm^*,G])^*$ for admissible cohomology
theories $E$, as we now explain.  Given such a cohomology theory we
need to define an associated ring $D'$.  In the special case of Morava
$E$-theory, this was defined in~\cite{hokura:ggc}; the details
necessary for the general case are given in~\cite{st:fsf}.  Associated
to $E$ we have a formal group $\GG$ over $\spf(E^0)$ and thus a
level-structure scheme $\Level(\Lm(m),\GG)$ with coordinate ring $D_m$
say.  These form a directed system in an obvious way and we define
$D'=\Q\ot\colim_mD_m$.  (This was called $L$ in~\cite{hokura:ggc} but
we have renamed it to avoid clashes of notation.)  If $\GG$ is the
universal deformation of its restriction to the special fibre (as is
the case with Morava $E$-theory) then $D'$ is the integral domain
obtained from $\Q\ot E^0$ by adjoining a full set of roots of
$[p^m](x)$ for all $m$.  For any $E$ one can show that $D'$ is a free
module of countable rank over $\Q\ot E^0$.

As mentioned earlier, the following theorem is merely a repackaging of
results of Hopkins, Kuhn and Ravenel~\cite{hokura:ggc}.
\begin{theorem}\label{thm-hkr}
 For any admissible ring spectrum $E$, there is a natural isomorphism
 of Frobenius algebras over $D'$
 \[ D'\ot_{E^0} E^0BG = D'\ot_\Q C([\Lm^*,G])^*. \]
 Moreover, this respects the constructions $u\mapsto(Lu)^*$ and
 $u\mapsto(Ru)^*$ for functors between groupoids.
\end{theorem}
\begin{proof}
 We first construct a map 
 \[ \tau\:D'\ot_{E^0}E^0BG\xra{}D'\ot_\Q C([\Lm^*,G])^* \]
 of $D'$-algebras.  By juggling various adjunctions we see that it
 suffices to construct, for each functor $\lm\:\Lm^*\xra{}G$, a map
 $\tau_\lm\:E^0BG\xra{}D'$ of $E^0$-algebras, such that
 $\tau_\lm=\tau_\mu$ when $\lm$ is isomorphic to $\mu$.  We know from
 our previous remarks that $\lm$ must factor through
 $\Lm(m)^*=\Lm^*/p^m$ for some $m$.  We thus get a map
 $E^0BG\xra{}E^0B\Lm(m)^*$, and we know from~\cite[Proposition
 2.9]{grst:vlc} that $E^0B\Lm(m)^*=\O_{\Hom(\Lm(m),\GG)}$, and $D_m$
 is a quotient of this ring, so we get the required map $\tau_\lm$ as
 the composite
 \[ E^0BG \xra{B\lm^*} E^0B\Lm(m)^* \xra{} D_m \xra{} D'. \]
 One checks easily that this is independent of the choice of $m$.
 Isomorphic functors $\lm,\mu$ give homotopic maps $B\Lm(m)^*\xra{}BG$
 and thus $\tau_\lm=\tau_\mu$ as required.  The resulting map $\tau$
 is easily seen to be natural for functors of groupoids and to convert
 equivalences to isomorphisms.  Both source and target of $\tau$
 convert disjoint unions to products.  Any finite groupoid is
 equivalent to a finite disjoint union of finite groups, so it
 suffices to check that $\tau$ is an isomorphism when $G$ is a group.
 This is just~\cite[Theorem B]{hokura:ggc}.

 To say that this isomorphism respects the construction
 $u\mapsto(Lu)^*$ is just to say that $\tau$ is a natural map, which
 is clear.  We also need to check that for any functor $u\:G\xra{}H$,
 the following diagram commutes:
 \begin{diag}
  \node{D'\ot_{E^0}E^0BG} \arrow{e,t}{\tau_G} \arrow{s,l}{(Ru)^*}
  \node{D'\ot_\Q C([\Lm^*,G])^*}              \arrow{s,r}{(Ru)^*}\\
  \node{D'\ot_{E^0}E^0BH} \arrow{e,b}{\tau_H}
  \node{D'\ot_\Q C([\Lm^*,H])^*}             
 \end{diag}
 We first make this more explicit.  The functor $u$ induces
 $Ru\:LH\xra{}LG$.  By applying $E^0(-)$ and noting that $E^0LK=E^0BK$
 we get a map $(Ru)^*\:E^0BG\xra{}E^0BH$.  After tensoring with $D'$
 we obtain the left hand vertical map in the above diagram.  On the
 other hand, $u$ also induces a functor
 $u_*\:[\Lm^*,G]\xra{}[\Lm^*,H]$ and thus a map
 $R(u_*)\:C[\Lm^*,H]\xra{}C[\Lm^*,G]$.  By dualising and tensoring
 with $D'$ we obtain the right hand vertical map.

 We first prove that the diagram commutes when $u$ is a
 quasi-covering.  This reduces easily to the case where $H$ is a group
 and $G$ is connected.  It is not hard to see that in this case $u$ is
 equivalent to the inclusion of a subgroup $G\leq H$ and
 $Ru\:LH\xra{}LG$ is just the $K(n)$-localisation of the transfer map
 $\Sgi BH_+\xra{}\Sgi BG_+$.  It follows
 from~\cite[Proposition 3.6.1]{hokura:ggc} that
 $\tau_\lm((Ru)^*x)=\sum\tau_{\lm^h}(x)$, where the sum runs over
 cosets $hG$ such that $\lm^h:=h^{-1}\lm h$ maps $\Lm^*$ into $G$.
 The right hand side can be rewritten as
 $|G|^{-1}\sum_h\tau_{\lm^h}(x)$, where the sum now runs over elements
 rather than conjugacy classes.  Fix a homomorphism
 $\mu\:\Lm^*\xra{}G$ that becomes conjugate to $\lm$ in $H$.  Then the
 number of $h$'s for which $\lm^h=\mu$ is the order of the group
 $Z_H(\lm)=\{h\in H\st \lm^h=\lm\}$, so 
 \[ \tau_\lm((Ru)^*x)=|G|^{-1}\sum_\mu|Z_H(\lm)|\tau_\mu(x). \]
 If we want to index this sum using conjugacy classes of $\mu$'s
 rather than the $\mu$'s themselves, we need an extra factor of
 $|G|/|Z_G(\mu)|$, the number of conjugates of $\mu$.  This gives
 \[ \tau_\lm((Ru)^*x)=
     \sum_{[\mu]}|Z_G(\mu)|^{-1}|Z_H(\lm)|\tau_\mu(x). 
 \]
 On the other hand, $Z_G(\mu)$ is just the automorphism group of $\mu$
 in the category $[\Lm^*,G]$, so the map
 $R(u_*)\:C[\Lm^*,H]\xra{}C[\Lm^*,G]$ is given by
 \[ R(u_*)[\lm]/|Z_H(\lm)| =
     \sum_{[\mu]\st u\mu\simeq\lm} [\mu]/|Z_G(\mu)|.
 \]
 The claim follows easily by comparing these formulae.
  
 We have an inner product on $D'\ot_{E^0}E^0BG$ obtained from the
 inner product $b_G$ on $LG$, and an inner product on
 $D'\ot_\Q C([\Lm^*,G])^*$ obtained from the standard inner product on
 $C(K)^*$ for any $K$.  By taking $u$ to be the diagonal functor
 $G\xra{}G\tm G$ in the above discussion, we see that our isomorphism
 $\tau$ converts the former inner product to the latter one.  Thus
 $\tau$ is compatible with taking adjoints and with the construction
 $u\mapsto(Lu)^*$, so it is compatible with the construction
 $u\mapsto(Ru)^*$ as well.
\end{proof}

We next reformulate Theorem~\ref{thm-hkr} in the spirit
of~\cite[Theorem 3.7]{grst:vlc}.  
\begin{definition}
 Given a finite groupoid $G$, we define a new groupoid $AG$ as
 follows.  The objects are pairs $(a,A)$, where $a\in G$ and $A$ is a
 finite Abelian $p$-subgroup of $G(a)$.  The morphisms from $(a,A)$ to
 $(b,B)$ are maps $g\:a\xra{}b$ in $G$ such that $B=gAg^{-1}$.  For
 any finite Abelian $p$-group $A$ we can define a ring
 $D'_A=\Q\ot\O_{\Level(A^*,\GG)}$ as in~\cite[Proposition 22]{st:fsf}.
 There is an evident way to make the assignment $(a,A)\mapsto D'_A$
 into a functor $AG^{\text{op}}\xra{}\text{Rings}$, and we define
 $TG=\invlim_{(a,A)\in AG}D'_A$.  If we write
 $\al'_{(a,A)}=|Z_{G(a)}(A)|$ then $\al'\in TG$.
\end{definition}

\begin{theorem}
 There is a natural isomorphism $\Q\ot E^0BG=TG$, and this is a
 finitely generated free module over $E^0$.  The element
 $\al=\mu\psi(1)\in E^0BG$ becomes $\al'\in TG$, so the resulting
 Frobenius form on $TG$ is just $\ep(x)=\tht(x/\al')$, where $\tht$ is
 the trace form.
\end{theorem}
\begin{proof}
 The isomorphism $\Q\ot E^0BG=TG$ can be proved either by reducing to
 the case of a group and quoting~\cite[Theorem 3.7]{grst:vlc}, or by
 taking the fixed-points of both sides in Theorem~\ref{thm-hkr} under
 the action of $\Aut(\Lm)$.  From the latter point of view, the term
 in $TG$ indexed by $(a,A)$ corresponds to the terms in $C[\Lm^*,G]^*$
 coming from homomorphisms $\Lm^*\xra{}G(a)$ with image $A$, so $\al'$
 becomes the function $\pi_0[\Lm^*,G]\xra{}L$ that sends $[\lm]$ to
 $|[\Lm^*,G](\lm)|$.  Proposition~\ref{prop-trace-loops} identifies
 this with $\al$, as required.
\end{proof}

\section{Warnings}
\label{sec-warnings}

We started this paper by considering the representation ring $R(G)$,
but unfortunately the analogy between our rings $E^0LG$ and $R(G)$
fails in a number of respects, even in the height one case.  In this
section we point out some possible pitfalls.

Let $E$ be the $p$-adic completion of the complex $K$-theory spectrum,
so $E$ is an admissible ring spectrum of height one.  Then $E^0BG$ is
the completion of $R(G)$ at $I+(p)$, where $I$ is the augmentation
ideal.  The ring $R(G)$ is a free Abelian group of rank equal to the
number of conjugacy classes, generated by the irreducible characters.
These are orthonormal, so the inner product on $R(G)$ is equivalent to
the standard diagonal, positive definite inner product on $\Z^h$.  It
also follows that $R(G)$ is a permutation module for the outer
automorphism group of $G$.

The ring $E^0BG$ is a free module over $\Zp$ of rank equal to the
number of conjugacy elements of elements of $p$-power order.  The
canonical map $R(G)\xra{}E^0BG$ does not preserve inner products.
There is no canonical set of generators for $E^0BG$, so there is no
reason for it to be a permutation module for $\Out(G)$.  In fact, Igor
Kriz has constructed examples of extensions $G\xra{}G'\xra{}C_p$ where
$G$ is good but $H^1(C_p;E^0BG)\neq 0$ and one can deduce that $E^0BG$
is not a permutation module in this case~\cite{kr:mkc}.  There
is also no reason to expect that $E^0BG$ has an orthonormal basis.

A related set of issues involves the comparison between the free loop
space of $BG$ (which is $B[\Z,G]$) and the space of maps from the
$p$-adically completed circle to $BG$ (which is $B[\Zp,G]$).  The
former enters into Proposition~\ref{prop-trace-loops}, and the latter
into Theorem~\ref{thm-hkr}.  The two spaces are not even $p$-adically
equivalent: if $G$ is a group and $T$ is a set of representatives for
the conjugacy classes of elements whose order is not a power of $p$
then $B[\Z,G]\simeq B[\Zp,G]\amalg\coprod_{g\in T}BZ_G(g)$, and each
term in the coproduct contributes at least a factor of $E^0$ in
$E^0B[\Z,G]$ even if $Z_G(g)$ is a $p'$-group.  Note
that~\cite[5.3.10]{hokura:ggc} is slightly inaccurate in this regard;
the proof given there really shows that
$\chi_nBG=\chi_{n-1}(B[\Zp,G])$, rather than $\chi_{n-1}(B[\Z,G])$.


\begin{thebibliography}{10}

\bibitem{ab:tdt}
L.~Abrams.
\newblock Two-dimensional topological quantum field theories and {F}robenius
  algebras.
\newblock {\em J. Knot Theory Ramifications}, 5(5):569--587, 1996.

\bibitem{an:fgr}
D.~W. Anderson.
\newblock Fibrations and geometric realizations.
\newblock {\em Bull. Amer. Math. Soc.}, 84(5):765--788, 1978.

\bibitem{bawe:ttt}
M.~Barr and C.~Wells.
\newblock {\em Toposes, Triples and Theories}, volume 278 of {\em Grundleheren
  der math. Wiss.}
\newblock Springer--Verlag, Berlin, 1985.

\bibitem{ba:ah}
H.~J. Baues.
\newblock {\em Algebraic Homotopy}, volume~15 of {\em Cambridge Studies in
  Advanced Mathematics}.
\newblock Cambridge University Press, 1989.

\bibitem{be:rcii}
D.~Benson.
\newblock {\em Representations and Cohomology {II}}, volume~31 of {\em
  Cambridge Studies in Advanced Mathematics}.
\newblock Cambridge University Press, 1991.

\bibitem{bo:lspe}
A.~K. Bousfield.
\newblock The localization of spectra with respect to homology.
\newblock {\em Topology}, 18:257--281, 1979.

\bibitem{ca:cc}
S.~Carnody.
\newblock {\em Cobordism categories}.
\newblock PhD thesis, Cambridge University, 1995.

\bibitem{dwsp:htm}
W.~G. Dwyer and J.~Spalinski.
\newblock Homotopy theories and model categories.
\newblock In I.~M. James, editor, {\em Handbook of Algebraic Topology}.
  Elsevier, Amsterdam, 1995.

\bibitem{ekmm:rma}
A.~D. Elmendorf, I.~Kriz, M.~A. Mandell, and J.~P. May.
\newblock {\em Rings, Modules and Algebras in Stable Homotopy Theory},
  volume~47 of {\em Amer. Math. Soc. Surveys and Monographs}.
\newblock American Mathematical Society, 1996.

\bibitem{frpi:cst}
R.~Fritsch and R.~A. Piccinini.
\newblock {\em Cellular structures in topology}, volume~19 of {\em Cambridge
  studies in advanced mathematics}.
\newblock Cambridge University Press, 1990.

\bibitem{frta:ant}
A.~Fr{\"o}lich and M.~Taylor.
\newblock {\em Algebraic number theory}.
\newblock Number~27 in Cambridge Studies in Advanced Mathematics. Cambridge
  University Press, 1991.

\bibitem{grma:gtc}
J.~P.~C. Greenlees and J.~P. May.
\newblock {\em Generalized {T}ate Cohomology}, volume 113 of {\em Memoirs of
  the American Mathematical Society}.
\newblock American Mathematical Society, 1995.

\bibitem{grst:vlc}
J.~P.~C. Greenlees and N.~P. Strickland.
\newblock Varieties and local cohomology for chromatic group cohomology rings.
\newblock 40 pp., To appear in Topology, 1996.

\bibitem{hokura:ggc}
M.~J. Hopkins, N.~J. Kuhn, and D.~C. Ravenel.
\newblock Generalised group characters and complex oriented cohomology
  theories.
\newblock Preprint (various editions, since the late 1980's).

\bibitem{hokura:mkc}
M.~J. Hopkins, N.~J. Kuhn, and D.~C. Ravenel.
\newblock Morava {K}-theories of classifying spaces and generalized characters
  for finite groups.
\newblock In J.~Aguad\'e, M.~Castellet, and F.~R. Cohen, editors, {\em
  Algebraic Topology: Homotopy and Group Cohomology}, volume 1509 of {\em
  Lecture Notes in Mathematics}, pages 186--209, New York, 1992.
  Springer-Verlag.

\bibitem{hopast:ash}
M.~Hovey, J.~H. Palmieri, and N.~P. Strickland.
\newblock Axiomatic stable homotopy theory.
\newblock {\em Mem. Amer. Math. Soc.}, 128(610):x+114, 1997.

\bibitem{host:mkl}
M.~Hovey and N.~P. Strickland.
\newblock Morava ${K}$-theories and localisation.
\newblock {\em Mem. Amer. Math. Soc.}, 139(666):104, 1999.

\bibitem{jost:gtci}
A.~Joyal and R.~Street.
\newblock The geometry of tensor calculus, {I}.
\newblock {\em Advances in Mathematics}, 88:55--112, 1991.

\bibitem{kr:mkc}
I.~Kriz.
\newblock Morava ${K}$-theory of classifying spaces: some calculations.
\newblock {\em Topology}, 36(6):1247--1273, 1997.

\bibitem{la:ccc}
P.~S. Landweber.
\newblock Complex cobordism of classifying spaces.
\newblock {\em Proceedings of the American Mathematical Society}, 27:175--179,
  1971.

\bibitem{lemast:esh}
L.~G. Lewis, J.~P. May, and M.~S. (with contributions~by Jim E.~McClure).
\newblock {\em Equivariant Stable Homotopy Theory}, volume 1213 of {\em Lecture
  Notes in Mathematics}.
\newblock Springer--Verlag, New York, 1986.

\bibitem{qu:haki}
D.~Quillen.
\newblock Higher algebraic ${K}$-theory. {I}.
\newblock In {\em Algebraic $K$-theory, I: Higher $K$-theories (Proc. Conf.,
  Battelle Memorial Inst., Seattle, Wash., 1972)}, pages 85--147. Lecture Notes
  in Math., Vol. 341, Berlin, 1973. Springer.

\bibitem{qu:ha}
D.~G. Quillen.
\newblock {\em Homotopical Algebra}, volume~43 of {\em Lecture Notes in
  Mathematics}.
\newblock Springer--Verlag, 1967.

\bibitem{qu:grk}
D.~G. Quillen.
\newblock The geometric realization of a {K}an fibration is a {S}erre
  fibration.
\newblock {\em Proc. Amer. Math. Soc.}, 19:1499--1500, 1968.

\bibitem{qu:epc}
D.~G. Quillen.
\newblock Elementary proofs of some results of cobordism theory using
  {S}teenrod operations.
\newblock {\em Advances in Mathematics}, 7:29--56, 1971.

\bibitem{ra:lrc}
D.~C. Ravenel.
\newblock Localization with respect to certain periodic homology theories.
\newblock {\em American Journal of Mathematics}, 106:351--414, 1984.

\bibitem{st:fsfg}
N.~P. Strickland.
\newblock Formal schemes and formal groups.
\newblock In J.~Meyer, J.~Morava, and W.~Wilson, editors, {\em
  Homotopy-invariant algebraic structures: in honor of J.M. Boardman},
  Contemporary Mathematics. American Mathematical Society, 1989.
\newblock 75 pp., To Appear.

\bibitem{st:fsf}
N.~P. Strickland.
\newblock Finite subgroups of formal groups.
\newblock {\em Journal of Pure and Applied Algebra}, 121:161--208, 1997.

\bibitem{st:pmm}
N.~P. Strickland.
\newblock Products on {$MU$}-modules.
\newblock Trans. Amer. Math. Soc., posted on March 1, 1999, PII:S
  0002-9947(99)02436-8 (31 pp., to appear in print), 1999.

\bibitem{sttu:rme}
N.~P. Strickland and P.~R. Turner.
\newblock Rational {M}orava {$E$}-theory and {$DS^0$}.
\newblock {\em Topology}, 36(1):137--151, 1997.

\bibitem{ti:csi}
U.~Tillmann.
\newblock The classifying space of the $(1+1)$-dimensional cobordism category.
\newblock {\em J. Reine Angew. Math.}, 479:67--75, 1996.

\bibitem{tu:dlo}
P.~R. Turner.
\newblock {D}yer-{L}ashof operations in the {H}opf ring for complex cobordism.
\newblock {\em Mathematical Proceedings of the Cambridge Philosophical
  Society}, 114:453--460, 1993.

\end{thebibliography}

\end{document}